\documentclass[11pt]{article}

\usepackage{graphicx}
\usepackage{amsmath}
\usepackage{amsthm}
\usepackage{amssymb}
\usepackage{enumerate}
\usepackage[breaklinks]{hyperref}
\usepackage{lscape}
\usepackage{multirow}
\usepackage{relsize}
\usepackage{setspace}
\usepackage[comma]{natbib}
\usepackage{fullpage}
\usepackage{caption}
\usepackage{authblk}
\usepackage{totcount}
\usepackage[titletoc]{appendix}
\usepackage{bigstrut}
\usepackage{color}
\usepackage{colortbl}
\usepackage[top    = 0.5in,
bottom = 1in,
left   = 1in,
right  = 1in]{geometry}
\newtotcounter{citnum} 
\def\oldbibitem{} \let\oldbibitem=\bibitem
\def\bibitem{\stepcounter{citnum}\oldbibitem}

\definecolor{Grey}{gray}{0.95}

\begin{document}
\onehalfspacing

\title{A New Simheuristic Approach for Stochastic Runway Scheduling}

\vspace{4mm}

\author{\large Rob Shone\footnote{Corresponding author}\\ \vspace{-4mm} \footnotesize{Centre for Transport and Logistics (CENTRAL), Lancaster University, Bailrigg, Lancaster, UK. LA1 4YW.\\  Email: r.shone@lancaster.ac.uk}\\\text{ }\\ \large Kevin Glazebrook\\ \footnotesize{Centre for Transport and Logistics (CENTRAL), Lancaster University, Bailrigg, Lancaster, UK.  LA1 4YW.\\  Email: k.glazebrook@lancaster.ac.uk}\\\text{ }\\ \large Konstantinos G. Zografos\\ \footnotesize{Centre for Transport and Logistics (CENTRAL), Lancaster University, Bailrigg, Lancaster, UK.  LA1 4YW.\\  Email: k.zografos@lancaster.ac.uk}}
\date{ }


\normalsize
\maketitle

\hrule
\text{ }\\

\noindent \Large \textbf{Abstract}\\
\normalsize

\noindent We consider a stochastic, dynamic runway scheduling problem involving aircraft landings on a single runway.  Sequencing decisions are made with knowledge of the estimated arrival times (ETAs) of all aircraft due to arrive at the airport, and these ETAs vary according to continuous-time stochastic processes.  Time separations between consecutive runway landings are modeled via sequence-dependent Erlang distributions and are affected by weather conditions, which also evolve continuously over time.  The resulting multi-stage optimization problem is intractable using exact methods and we propose a novel simheuristic approach, based on the application of methods analogous to variable neighborhood search (VNS) in a high-dimensional stochastic environment.  Our model is calibrated using flight tracking data for over 98,000 arrivals at Heathrow Airport.  Results from numerical experiments indicate that our proposed simheuristic algorithm outperforms an alternative based on deterministic forecasts under a wide range of parameter values, with the largest benefits being seen when the underlying stochastic processes become more volatile and also when the on-time requirements of individual flights are given greater weight in the objective function. \\

\noindent \textit{Keywords: }Simulation optimization, stochastic processes, runway scheduling, aviation\\

\hrule

\newpage

\text{ }\\

\section{Introduction}\label{sec:intro}

\text{ }\indent Imbalances between demand and capacity at the world's busiest airports continue to pose problems for schedule coordinators, air traffic controllers, airspace users and other stakeholders.  By the end of September 2022, daily flight numbers in Europe had recovered to $88\%$ of their pre-pandemic levels and major hubs such as London Heathrow, Paris Charles de Gaulle and Amsterdam Schiphol were again processing more than 1000 runway movements (i.e. take-offs or landings) per day on average (\cite{Eurocontrol2022a}).  Flight delays, however, are also increasing in frequency following the record lows seen during 2020 (\cite{Eurocontrol2022b}).  The question of how to mitigate air traffic congestion while satisfying the ever-increasing demand for air transport services has been examined extensively at the strategic, tactical and operational levels (\cite{Zografos2017,Jacquillat2018,Cavusoglu2021}).

\text{ }\indent Air traffic delays can be mitigated by improving the efficiency of airport operations, including runway sequencing.  In this paper we consider a stochastic, dynamic runway scheduling problem.  The focus is on tactical decision-making\footnote{The term `operational decision-making' may be preferred by some authors.}, in the sense that we aim to optimize the real-time decisions made by air traffic controllers in response to the latest available information on a particular day of operations, with the schedule of aircraft take-offs and landings having been determined in advance.  Our problem formulation draws upon different areas of the literature that have evolved along quite separate lines in the last few decades.  On one hand, we model the stochastic nature of airport runway operations using a queueing theory approach, whereby the times that aircraft arrive at the runway threshold and complete `service' (i.e. usage of the runway) are subject to uncertainty.  Many previous studies have used stochastic queueing formulations to model operational delays at airports (\cite{Koopman1972,Stamatopoulos2004,Stolletz2008,Hansen2009,Pyrgiotis2016,Jacquillat2017}).  The decision-making aspect of our problem relates to the sequencing of aircraft landing at a single runway during a congested period.  In this respect, we build upon the literature on aircraft sequencing problems, which began with the formulation of static, deterministic optimization problems (\cite{Psaraftis1978}) but has more recently expanded into the area of stochastic optimization (\cite{Solveling2011,Heidt2016,Liu2018,Solak2018,Khassiba2020}).

\text{ }\indent In queueing models of airport runway operations, a common approach is to use an estimate for the airport's capacity in order to derive `service rates' for the queues.  The capacity of an airport or single runway may be defined as the expected number of runway movements per unit time that can be operated under conditions of continuous demand (\cite{Neufville2013}).  This definition implies that an airport's capacity actually varies with time, as it depends on various controllable and non-controllable factors, including runway configurations and weather conditions.  \cite{Gilbo1993} introduced the concept of a `capacity envelope' to represent the set of feasible pairs of service rates for arrivals and departures at a single airport, and subsequently this approach has been used to formulate stochastic, dynamic optimization problems based on the control of service rates at discrete time epochs (\cite{Jacquillat2015,Jacquillat2017,Shone2019}).  The use of a capacity envelope to select service rates for aircraft queues may be seen as somewhat macroscopic in nature, as it does not explicitly allow for fine-grain aircraft sequencing and the time savings that air traffic controllers can achieve by taking into account separation requirements between different types of aircraft; instead, it assumes that different possible traffic mixes and other considerations can be implicitly accounted for by configuring service time variances (\cite{Shone2021}).  The model proposed in this paper assumes that an aircraft's service time follows a random distribution which depends on its own weight class and also that of its immediate predecessor in the runway queue.  Thus, we allow for the effects of weight classes on separation times, while preserving the stochastic modeling of runway service times.

\text{ }\indent Two-stage stochastic optimization has been popularized in recent years as a means of incorporating uncertainty into runway scheduling problems.  \cite{Solveling2011} introduced a two-stage model in which a sequence of aircraft weight classes is determined in the first stage, and specific flights are assigned to positions in the sequence (subject to weight class compatibility) in the second stage.  \cite{Solak2018} subsequently enhanced this model by introducing costs based on exact timings of runway operations.  \cite{Liu2020} used a similar formulation, but allowed for limited or ambiguous information about the model parameters.  \cite{Khassiba2020} considered the optimization of a sequence of aircraft arriving at an initial approach fix (IAF) and used chance constraints to mitigate the risk of separation time violations.  A standard assumption in two-stage optimization problems is that the first-stage decisions are made under uncertainty, but the second stage decisions are made after all uncertainty has been realized.  As such, uncertainty is `on/off' in nature.  The two-stage methodology therefore appears best-suited to problem instances involving relatively small numbers of aircraft and short time horizons, in which it makes sense for the uncertainties associated with all aircraft (related to their arrival times at an IAF, for example) to be realized at the same point in time.  One can, however, use a `rolling horizon' approach in order to apply the two-stage method to larger problem instances.

\text{ }\indent Two-stage optimization methods clearly have useful applications in runway scheduling problems, and offer the additional advantage of being solvable using well-established optimization techniques.  In this paper, however, we take a different approach in order to model the problem faced by air traffic controllers who need to make decisions in information-rich and fast-changing dynamic environments.  We consider multi-stage problems involving hundreds of aircraft, in which hundreds of thousands of decisions are made over an operational period of a few hours.  Furthermore, we use continuous-time stochastic processes to model the evolution of uncertainty over time, with each decision being made under the latest set of information available.  The high-dimensional nature of our problem precludes the use of any exact solution methods.  Instead, we opt for a simulation-based approach in which many possible runway sequences are compared using randomly-sampled system trajectories.  We continuously update `fitness estimates' for the best-performing sequences found so far and aim to discover new ones by adapting the search neighborhood according to the latest system state.  Our solution approach may be recognized as a form of \emph{simheuristic}; see \cite{Juan2015} for some useful background on this fast-developing research area.

\text{ }\indent The main contributions of our paper are as follows:\\

\begin{itemize}
\item We provide a new formulation for the multi-objective stochastic runway scheduling problem which includes three different types of dynamic uncertainty: (i) estimated times of arrival (ETAs) for aircraft evolve according to continuous-time stochastic processes; (ii) sequence-dependent aircraft separation times follow Erlang distributions; (iii) expected changepoints in weather conditions also vary according to continuous-time stochastic processes.
\item We propose a novel solution methodology for this problem, based on the application of simheuristic search techniques in a stochastic and rapidly-changing environment.
\item We configure our model using data from over 98,000 flights landing at Heathrow Airport in 2018 and 2019 and show, using numerical experiments, that the solutions given by our simheuristic approach consistently outperform those given by alternative heuristics with respect to a dual objective function based on schedule punctuality and air-holding times.\\
\end{itemize}

\text{ }\indent Our model formulation is provided in Section \ref{sec:formulation}.  Details of our heuristic approaches are presented in Section \ref{sec:methodology}.  The use of flight tracking data to configure our model parameters is described in Section \ref{sec:data}, and we provide details of results from our computational experiments in Section \ref{sec:numerical}.  Our concluding remarks are given in Section \ref{sec:conclusions}.\\

\section{Model formulation}\label{sec:formulation}

\text{ }\indent In this paper we restrict attention to arriving flights (landings) at a single airport.  We assume that all arrivals take place on the same runway, and that this runway is not used or affected by the stream of outbound traffic (departures).  This situation is quite common at airports which operate two parallel runways in `segregated mode', with Heathrow Airport in London being a pertinent example.  
Let $\mathcal F$ be the set of arriving flights scheduled to use the runway during a particular time interval, denoted by $\mathcal T:=[0,T]$.  We assume here that $\mathcal T$ is not longer than one day, so we are optimizing decisions over a period of hours rather than days or weeks.  For each flight $i\in\mathcal F$ we associate a scheduled arrival time $a_i\in\mathcal T$ at the destination airport and a scheduled departure time $d_i\in(-\infty,a_i)$ from its origin airport.  For clarity, we emphasize that the destination airport is the same for all flights in $\mathcal F$, but the origin airport is flight-specific.  We also use $w_i$ to denote the weight class of aircraft $i\in\mathcal F$, belonging to a set of weight classes $\mathcal W$, and we define $g_i\in [0,1]$ as a relative cost parameter associated with delays to flight $i$.

\text{ }\indent The actual times that flights land on the runway are affected not only by sequencing and scheduling decisions, but also by the uncertainty affecting (i) their departure and flight times, (ii) landing time separations with preceding aircraft and (iii) weather conditions at the destination airport.  In Sections \ref{sec21}-\ref{sec23} we explain how these different sources of uncertainty are modeled, and in Section \ref{sec24} we present our decision-making framework and objective function.\\

\subsection{Unconstrained landing times}\label{sec21}

\text{ }\indent The first source of uncertainty in our model is related to the earliest time that a flight would be able to land in the absence of congestion effects or adverse weather conditions at the destination airport.  Following previous studies (\cite{Bennell2017,Khassiba2020}), we refer to this as the \emph{unconstrained landing time} and denote it by $A_i$ for flight $i\in\mathcal F$.  Noting that $A_i$ may depend on many unpredictable factors and control interventions at various different stages of flight $i$'s progress (including the pre-departure stage), we propose to make a distinction between \emph{pre-tactical uncertainty} and \emph{tactical uncertainty} and write
\begin{equation}A_i=a_i+\Delta_i^{\text{pre}}+\Delta_i^{\text{tac}},\label{eq1}\end{equation}
where $\Delta_i^{\text{pre}}$ and $\Delta_i^{\text{tac}}$ are pre-tactical and tactical delays (possibly negative-valued), respectively.  Here, `pre-tactical' delays are those which can already be foreseen well in advance of a flight's departure, but would not have been known when the airport arrival schedule was originally produced (typically several months in advance of operations).  Examples of delays which may fall into the `pre-tactical' category are those caused by airline crew unavailability, local airspace restrictions or global upper wind conditions.  On the other hand, `tactical' delays are those which evolve dynamically during actual operations.  These might include take-off delays caused by taxiway congestion or enroute delays caused by the need to avoid potential air traffic conflicts.  

\text{ }\indent The pre-tactical delay $\Delta_i^{\text{pre}}$ is a semi-bounded continuous random variable in our model, while $\Delta_i^{\text{tac}}$ is also stochastic but additionally depends on sequencing decisions.  In our numerical experiments later, we rely on gamma distributions for modeling the pre-tactical delays, with 
flight-specific parameters estimated from historical data; further details are given in Section \ref{sec:data}.  For each $i\in\mathcal F$, we assume that $\Delta_i^{\text{pre}}$ is already `realized' in advance of flight $i$'s scheduled departure time and remains constant throughout the remainder of $\mathcal T$.  In practice, due to the way that our solution algorithms work, there is no loss of generality in assuming that $\Delta_i^{\text{pre}}$ is realized at the beginning of $\mathcal T$; see Section \ref{sec:methodology} for further details.  On the other hand, $\Delta_i^{\text{tac}}$ is not known until flight $i$ actually lands, and until then we are only able to predict its value.  

\begin{figure}[htbp]
    \begin{center}
        \includegraphics[width=16.5cm]{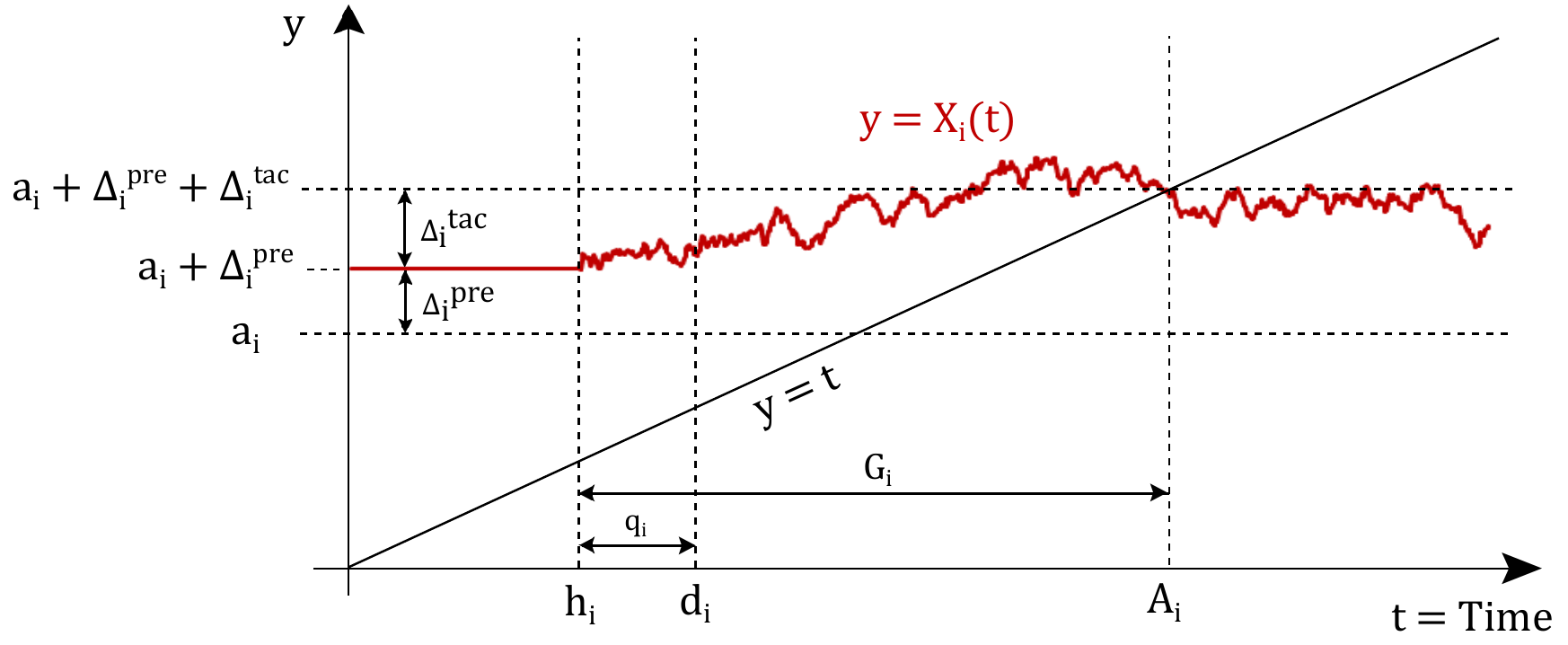}
    \end{center}
    \caption{The unconstrained landing time $A_i$ is obtained as the earliest point at which the trajectory $y=X_i(t)$, which behaves randomly for $t>h_i$, intersects the straight line $y=t$. The pre-tactical and tactical delays, $\Delta_i^{\text{pre}}$ and $\Delta_i^{\text{tac}}$, are shown as distances on the vertical scale.}\label{fig1}
\end{figure}

\text{ }\indent To model tactical uncertainty we use $X_i(t)$ to denote the \emph{estimated time of arrival} (ETA) for flight $i$ at time $t\in\mathcal T$, which depends on the latest information available at $t$.  The ETA $X_i(t)$ may vary considerably in short intervals of time and is generally not a monotonic function of $t$.  We assume $X_i(t)$ remains equal to $a_i+\Delta_i^{\text{pre}}$ (an adjusted ETA after realization of the pre-tactical delay $\Delta_i^{\text{pre}}$) until $t$ is within a certain proximity $q_i$ of the scheduled departure time $d_i$, at which point tactical uncertainty (which may include some pre-departure uncertainty associated with taxi-out times, for example) begins to take effect and $X_i(t)$ varies according to a Brownian motion (BM) process.  For convenience, let $h_i=d_i-q_i$.  Then we define
\begin{equation}X_i(t)=\begin{cases}a_i+\Delta_i^{\text{pre}},&\text{ if }t\leq h_i,\\
a_i+\Delta_i^{\text{pre}}+B_i(t)+\rho_i(t),&\text{ if }t>h_i,\end{cases}\label{BM_eq0}\end{equation}
where $B_i(t)\sim \text{N}(0,\;\sigma_i^2(t-h_i))$ and $\rho_i(t)\geq 0$ denotes any additional `airborne holding time' incurred by flight $i$ as part of the tactical decision-making process, i.e. the aircraft sequencing.  We elaborate further on this holding time in Section \ref{sec24}, but in the remainder of this subsection we assume (for ease of exposition) that $\rho_i(t)\equiv 0$.  The unconstrained landing time $A_i$ is given by
\begin{equation}A_i=\min\Big\{t>h_i\;:\;X_i(t)\leq t\Big\}.\label{BM_eq1}\end{equation}
\text{ }\indent That is, $A_i$ is the earliest point at which $X_i(t)$ is exceeded by the current time.  Figure \ref{fig1} shows how $A_i$ and $\Delta_i^{\text{tac}}$ are determined by the BM trajectory for a particular flight $i$, given some parameters $a_i,d_i,q_i$ and a fixed realization of $\Delta_i^{\text{pre}}$.  Note that we have the logical property that as $t$ increases, $X_i(t)$ becomes an increasingly accurate forecast of $A_i$; in other words, the unconstrained landing time becomes more predictable with time.

\text{ }\indent From (\ref{BM_eq1}), it follows that
\begin{align}A_i&\stackrel{\text{dist}}{=}\min\Big\{t>h_i\;:\;B_i(t)+t\geq a_i+\Delta_i^{\text{pre}}\Big\}\label{BM_eq2}\\
&=h_i+\min\Big\{u>0\;:\;B_i(u+h_i)+u\geq a_i+\Delta_i^{\text{pre}}-h_i\Big\}\label{BM_eq3}\\
&\stackrel{\text{def}}{=}h_i+G_i\label{BM_eq4},\end{align}
where we have used the fact that $-B_i(t)$ has the same distribution as $B_i(t)$ in (\ref{BM_eq2}), and the change of variables $u=t-h_i$ is used in (\ref{BM_eq3}).  Noting that $B_i(u+h_i)\sim\text{N}(0,\;\sigma_i^2 u)$, we see that $A_i$ can be interpreted as the `first hitting time' for a BM process with a positive linear drift (represented by the additive term $u$ in (\ref{BM_eq3})).  Using known theoretical results for BM processes with drift (see \cite{Cox1960} or \cite{Folks1978}), it follows that the conditional density function of the random variable $G_i$ defined in (\ref{BM_eq4}), given some realization of $\Delta_i^{\text{pre}}$, is the inverse Gaussian function

\begin{equation}f_{G_i}(u\;|\;\Delta_i^{\text{pre}}=\delta)=\frac{a_i+\delta-h_i}{\sigma_i\sqrt{2\pi u^3}}\exp\left(-\frac{(a_i+\delta-h_i-u)^2}{2\sigma_i^2 u}\right)\;\;\;\;(u>0),\label{BM_eq5}\end{equation}

\text{ }\\
which has mean $a_i+\delta-h_i$ and variance $\sigma_i^2(a_i+\delta-h_i)$.  The tactical delay $\Delta_i^{\text{tac}}$ can be expressed as
$$\Delta_i^{\text{tac}}=G_i-(a_i+\Delta_i^{\text{pre}}-h_i),$$
which can be recognized from Figure \ref{fig1} as the difference between the duration of the Brownian motion, $A_i-h_i$, and the duration that would occur if $X_i(t)$ remained at its initial level $a_i+\Delta_i^{\text{pre}}$ throughout $\mathcal T$.  Naturally, this has a mean of zero.  We also note, using (\ref{eq1}), that 
\begin{equation}\mathbb{E}[A_i]=a_i+\mathbb{E}[\Delta_i^{\text{pre}}]\label{BM_eq6a}\end{equation}
and, using the law of total variance,
\begin{align}\text{Var}(A_i)&=\text{Var}(G_i)\nonumber\\
&=\text{Var}(\mathbb{E}[G_i\;|\;\Delta_i^{\text{pre}}])+\mathbb{E}[\text{Var}(G_i\;|\;\Delta_i^{\text{pre}})]\nonumber\\
&=\text{Var}(\Delta_i^{\text{pre}})+\sigma_i^2\left(a_i+\mathbb{E}[\Delta_i^{\text{pre}}]-h_i\right).\label{BM_eq6b}\end{align}

\text{ }\indent Our use of BM trajectories to model the progress of flights towards their destination is not without precedent in the literature.  In problems based on aircraft conflict detection and resolution, for example, several authors have employed stochastic differential equations to represent the state evolution of an aircraft over time.  In these equations, the BM component is typically used to represent a stochastic wind disturbance; see, for example, \cite{Blom2002,Hu2005,Liu2014}.  Our model, on the other hand, uses BM trajectories in a rather more general sense to model variations in flight ETAs over time.

\text{ }\indent While it is possible to obtain data on the timeliness of aircraft arrivals and departures at airports around the world, it is not easy to find reliable information about how an arrival's ETA varies dynamically during (and prior to) its flight.  Hence, parameters such as $\sigma_i$ in our model are difficult to estimate accurately.  In our numerical experiments in Section \ref{sec:numerical}, we consider a range of possible values of $\sigma_i$, but for each value of $\sigma_i$ we adjust the parameters of the pre-tactical delay distribution in such a way that the value of $\text{Var}(A_i)$ given by (\ref{BM_eq6b}) corresponds as closely as possible to an estimate of $\text{Var}(A_i)$ obtained from historical data.  In other words, we keep $\text{Var}(A_i)$ close to a fixed, data-calibrated value, but experiment with different `weightings' for the pre-tactical and tactical components of the variance in (\ref{BM_eq6b}).\\

\subsection{Landing separation times}\label{sec22}

\text{ }\indent In our model, the actual time that flight $i\in\mathcal F$ lands on the runway depends not only on the unconstrained landing time $A_i$ but also on the traffic congestion and weather conditions at the destination airport.  In reality, pairs of aircraft landing consecutively on the same runway are required to maintain certain time separations, and these separations depend on the weight classes of the aircraft involved; for example, the gap between a leading aircraft and a following aircraft tends to be longer if the leader is in the `heavy' class, due to the extra wake turbulence generated (\cite{Newell1979}).  Many previous studies have incorporated class-dependent separation requirements in mathematical formulations of aircraft sequencing and scheduling problems, although these problems are often of a deterministic nature (\cite{Dear1976,Psaraftis1978,Beasley2000,Beasley2004,Bennell2017}).

\text{ }\indent Congestion in airport terminal areas naturally causes `queues' to form as aircraft await their turn to use the runway.  
There is a long-established tradition of using time-dependent queueing models such as $M(t)/E_k(t)/1$ to model on-time performance at airports (\cite{Kivestu1976,Malone1995,Stamatopoulos2004,Hansen2009,Pyrgiotis2016,Jacquillat2017,Shone2019}).  Here, $E_k$ denotes an Erlang-$k$ distribution, which is a form of gamma distribution.  This type of distribution is often favored in queueing models of air traffic due to its convenience and configurability; indeed, the parameters of such distributions can often be adjusted in order to ensure close resemblance to empirical distributions that might be observed in practice (\cite{Gupta2010}).

\text{ }\indent To the best of our knowledge, there is no precedent in the literature for the incorporation of class-dependent, Erlang-distributed separation times in an aircraft sequencing or runway scheduling problem.  We propose that such an approach makes sense, given that both class-dependent separation times and Erlang queue service times have (separately) been common features in previous air traffic models.  Let us define $L_i$ as the actual time that flight $i\in\mathcal F$ touches down on the runway (it is implied that $L_i\geq A_i$).  Suppose that flight $j\in\mathcal F$ immediately follows $i$ in the landing sequence, and let $e_{ij}$ denote the recommended time separation between a leading aircraft of class $w_i\in\mathcal W$ and a following aircraft of type $w_j\in\mathcal W$.  Then the random variable $M_{ij}$, interpreted as the time between the landings of $i$ and $j$ if these landings occur during a `congested period', is assumed to be Erlang-distributed with the probability density function
$$f_{M_{ij}}(t)=\frac{k\mu_{ij}(k\mu_{ij} t)^{k-1}\text{exp}(-k\mu_{ij} t)}{(k-1)!}\;\;\;\;\;(t>0),$$
where $\mu_{ij}=e_{ij}^{-1}$ and $k$ is an integer-valued shape parameter.  The mean of this distribution is $e_{ij}$ and the variance is $e_{ij}^2/k$; hence, large values of $k$ result in time separations that conform closely to the recommended values.  For clarity, we emphasize that our model allows the possibility of actual time separations being smaller than their recommended values (due to random variation), although large $k$ values would result in only small deviations.  In our numerical experiments in Section \ref{sec:numerical}, we experiment with different values of $k$ but choose $e_{ij}$ according to actual air traffic regulations.

\text{ }\indent The actual landing time of flight $j$ is then given by

$$L_j=\begin{cases}A_j,&\text{ if }j\text{ is the first plane to land during }\mathcal T,\\
\max\{A_j,\;L_i+M_{ij}\},&\text{ if }j\text{ is preceded by }i\text{ in the landing sequence.}\end{cases}$$
\text{ }\\
\text{ }\indent That is, flight $j$ cannot land before its unconstrained landing time $A_j$, but also must be appropriately separated from the preceding aircraft $i$.  If flight $j$ happens to arrive in the terminal area much later than flight $i$'s landing time (i.e. $A_j>>L_i$) then the random separation $M_{ij}$ essentially becomes irrelevant.  This is why we refer to $M_{ij}$ as the time between landings during a `congested period'.

\text{ }\indent In fact, the required separation $e_{ij}$ between two consecutively-arriving flights $i$ and $j$ may also have a time dependence in our model due to the effects of weather conditions in the terminal area at the time of flight $j$'s final approach.  We have thus far written $e_{ij}$ (without any time dependence) in order to simplify the notation, but this should be modified if weather variations are to be included.  Further details are given in Sections \ref{sec23} and \ref{sec24}.\\

\subsection{Weather conditions in the terminal area}\label{sec23}

\text{ }\indent Poor weather conditions at an airport can reduce its operational capacity by enforcing longer time separations between consecutive take-offs and landings.  The effects of weather conditions on airport capacity have been included in many different ways in previous research studies.  Notably, much of the literature on airport ground delay programs and air traffic flow management assumes probabilistic capacity profiles for airports and/or air sectors, with decisions being made under uncertainty as a result (\cite{Odoni1987,Richetta1993,Ball2001,Vossen2011,Corolli2017,Estes2020}).  Also, in problems related to the pre-tactical or dynamic allocation of airport capacity between arrivals and departures, a popular approach is to construct an airport `capacity envelope' (as mentioned in Section \ref{sec:intro}) whose shape depends on many operational factors, including runway configurations and prevailing weather conditions (\cite{Gilbo1993,Simaiakis2013,Jacquillat2015}).

\text{ }\indent In our model, given that we wish to consider aircraft sequencing decisions during a single day, it makes sense to assume that any periods of bad weather (and their timings) can be predicted with a high level of precision.  However, there may still be some uncertainty in the forecasts.  We consider two possible cases: (i) conditions are certain to remain `fine' throughout $\mathcal T$; (ii) conditions will be fine except for a single period of bad weather, denoted $\mathcal U\subset \mathcal T$, during which longer time separations are required.  (To use the correct terminology, we should note that `fine weather' implies visual meteorological conditions (VMC), whereas bad weather implies instrumental meteorological conditions (IMC); see \cite{Jacquillat2015} for further details.)  In the first case, the time separations $e_{ij}$ referred to in Section \ref{sec22} do not require any time dependence.  In the second case, however, we use $R_j$ to denote the time that flight $j\in\mathcal F$ begins the final stage of its journey to the runway and define the required separation between $i,j\in\mathcal F$ by
\begin{equation}E_{ij}(R_j)=\begin{cases}e_{ij},&\text{ if }R_j\in\mathcal T\setminus \mathcal U,\\
\phi e_{ij},&\text{ if }R_j\in\mathcal U,\end{cases}\label{Eij_eq}\end{equation}
where $\phi>1$ and $e_{ij}$ is interpreted as the `fine weather' value.  Then, similarly to before, $M_{ij}$ is Erlang-distributed with parameters $\mu_{ij}(R_j)=1/E_{ij}(R_j)$ and $k\in\mathbb{N}$.  

\begin{figure}[htbp]
    \begin{center}
        \includegraphics[width=14.5cm]{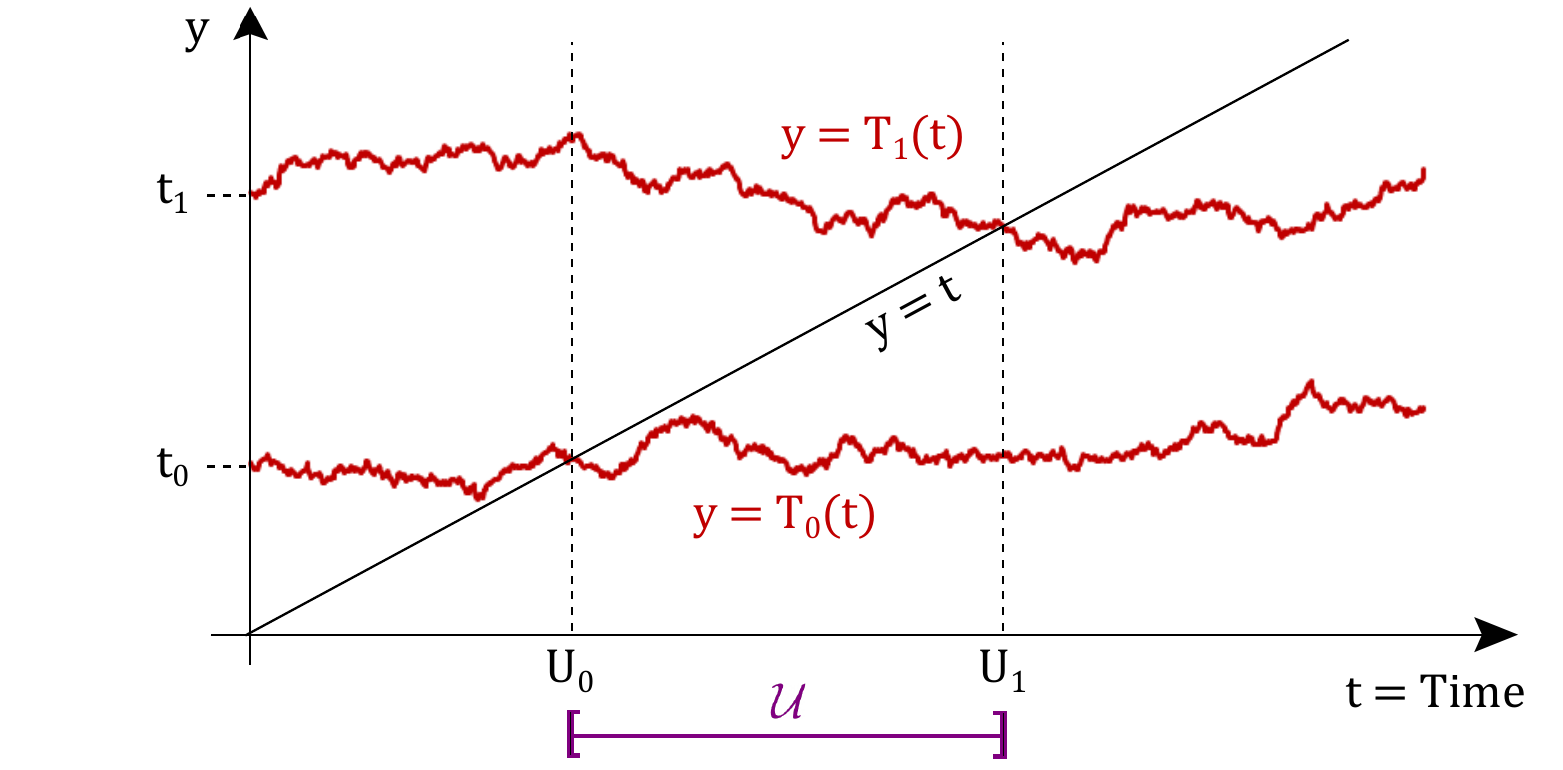}
    \end{center}
    \caption{The interval of bad weather, $\mathcal U$, is bounded between the intersection points $U_0=\min\{t>0\;|\;T_0(t)\leq t\}$ and $U_1=\min\{t>0\;|\;T_1(t)\leq t\}$.}\label{weather_fig}
\end{figure}

\text{ }\indent The beginning and ending times of $\mathcal U$ are subject to dynamic uncertainty.  Suppose that, at the beginning of $\mathcal T$, we expect that $\mathcal U$ will begin at time $t_0$ and end at some later time $t_1$ ($t_0,t_1\in\mathcal T$).  However, during $\mathcal T$ we continuously revise these estimates according to the latest forecast.  Let $T_0(t)$ and $T_1(t)$ be defined by Brownian motion trajectories as follows:
\begin{align}T_0(t)\sim\text{N} (t_0, \nu^2 t),\;\;\;\;\;\;T_1(t)\sim\text{N} (t_1,\nu^2 t),\nonumber\end{align}
where $\nu>0$ is a variance parameter.  We will assume independence between $T_0(t)$ and $T_1(t)$ for simplicity, although a dependence structure could be incorporated if needed.  Next, we define
\begin{align}U_0=\min\{t>0\;|\;T_0(t)\leq t\},\;\;\;\;\;\;U_1=\min\{t>0\;|\;T_1(t)\leq t\},\nonumber\end{align}
and the interval $\mathcal U$ is then given by
$$\mathcal U=\begin{cases}[U_0,\;U_1],&\text{ if }U_0<U_1,\nonumber\\
\emptyset,&\text{ otherwise.}\end{cases}$$
\text{ }\indent At any time $t$ satisfying $t<\min\{U_0,\;U_1\}$, our prediction is that $\mathcal U$ will begin at $T_0(t)$ and end at $T_1(t)$, assuming that $T_0(t)<T_1(t)$.  If $T_0(t)\geq T_1(t)$, then the prediction is that no period of bad weather will occur.  Figure \ref{weather_fig} provides an illustration of our approach.\\

\subsection{Decision-making framework and objective function}\label{sec24}

\text{ }\indent In this subsection we describe how our stochastic runway scheduling problem can be formulated as a stochastic, dynamic optimization problem by providing details of the state space, action space, cost mechanism and objective function.  First, however, we must provide some extra information about the statuses of flights in $\mathcal F$ and how these change during $\mathcal T$.

\text{ }\indent As explained in Section \ref{sec21}, we assume that the pre-tactical delays $\Delta_i^{\text{pre}}$ are known for all $i\in\mathcal F$ at the beginning of $\mathcal T=[0,T]$.  We scale the units of time in such a way that $0<a_i+\Delta_i^{\text{pre}}<T$ for all $i\in\mathcal F$; that is, all flights are expected to land during $\mathcal T$ following the realization of pre-tactical uncertainty.  It is possible that $h_i<0$ for some $i\in\mathcal F$, indicating that $X_i(t)$ is already varying according to Brownian motion and flight $i$ may (or may not) be airborne at $t=0$.  For any such flights, we generate an initial ETA $X_i(0)$ by sampling from a Normal distribution with mean $a_i+\Delta_i^{\text{pre}}$ and variance $\sigma_i^2 |h_i|$ in order to be consistent with (\ref{BM_eq0}).

\text{ }\indent For flight $i\in\mathcal F$, let $Q_i$ be defined by
$$Q_i=\min\Big\{t>0\;|\;X_i(t)-\tau\leq t\Big\},$$
where $\tau>0$ is a fixed parameter that defines the stage of an aircraft's journey at which it becomes eligible for tactical sequencing; in other words, sequencing decisions are made for aircraft that are within $\tau$ time units of their expected unconstrained landing time.  In practice, $\tau=30$ minutes might be an appropriate value (\cite{Bennell2017,Khassiba2020}).  At time $Q_i$, flight $i$ is said to enter a `pool' of aircraft waiting to be sequenced, i.e. given a position in the runway landing order.  The `pool' does not represent any particular physical location or region; instead, it merely represents a collection of aircraft that are expected to be able to land within $\tau$ time units.  In particular, it should not be confused with a holding stack, as planes in a holding stack must generally exit from the stack in a fixed order, whereas the `pool' in our model represents the last stage of an aircraft's journey at which its sequence position remains undetermined.

\text{ }\indent The sequencing decision for flight $i$ does not necessarily need to be made at time $Q_i$.  Instead, flight $i$ can be retained in the pool until some later time $R_i>Q_i$, at which point it is `released' and proceeds with the final stage of its journey.  During the time spent in the pool, the flight's progress is effectively `paused' and its ETA increases linearly, i.e. $X_i(t)=X_i(Q_i)+t-Q_i=t+\tau$ for $t\in[Q_i,R_i]$, indicating that it is not making further progress towards the runway during this time.  It may be performing path-stretching maneuvers, for example, in order to elongate its journey (\cite{Newell1979,Montlaur2017}).  
We define $\rho_i(t)=\min\{t,R_i\}-Q_i$ for $t>Q_i$ as the `holding time' incurred up to time $t$, satisfying equation (\ref{BM_eq0}).  It will also be convenient to define $\rho_i:=\lim_{t\rightarrow\infty}\rho_i(t)=R_i-Q_i$ as the total amount of time that flight $i$ spends in the pool before being released.

\begin{figure}[htbp]
    \begin{center}
        \includegraphics[width=15.5cm]{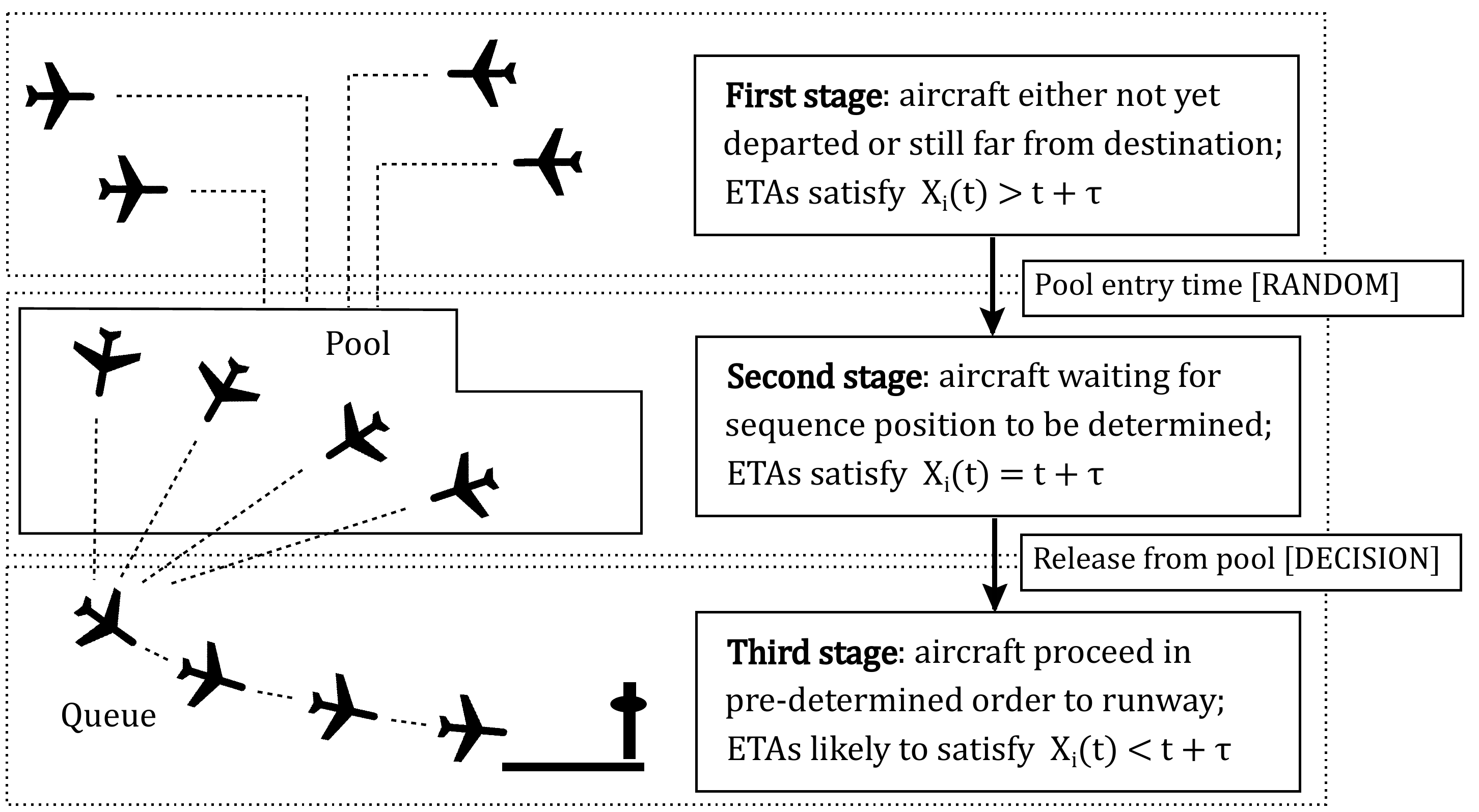}
    \end{center}
    \caption{Aircraft at different stages of their journeys at an arbitrary point in time $t\in\mathcal T$.}\label{sec24_fig}
\end{figure}

\text{ }\indent At time $R_i$, flight $i$ enters a `queue' of aircraft waiting to use the runway.  The queue is strictly first-come-first served, so that aircraft are required to land in the same order that they are released from the pool.  During periods of heavy congestion, the queue is likely to include planes circling in a holding stack, but planes on their final descent towards the runway are also notionally part of the `queue' in our model.  We note that
\begin{equation}A_i=\min\Big\{t>R_i\;|\;t\geq R_i+\tau+\tilde{B}_i(t)\Big\},\label{new_Ai_eq}\end{equation}
where $\tilde{B}_i(t)\sim\text{N}(0,\sigma_i^2(t-R_i))$.  This ensures consistency with the details in Section \ref{sec21}.  The actual landing time $L_i$ also depends on the required time separation between flight $i$ and its predecessor in the queue, as explained in Section \ref{sec22}, and this separation depends on the weather conditions at time $R_i$ as described in Section \ref{sec23}.  Figure \ref{sec24_fig} illustrates our model by showing aircraft at different stages of their journeys at some arbitrary point in time $t\in\mathcal T$.

\text{ }\indent We note that, due to the variable nature of the BM trajectories, it is possible for a flight's expected remaining flight time to exceed $\tau$ even if it has already entered the pool at an earlier time and has also been released from the pool; that is, we may have $X_i(t)>t+\tau$ for some $i\in\mathcal F$ even if $t>R_i$.  In this situation, we assume that flight $i$'s position in the queue (and, therefore, its position in the landing sequence) remains unaffected.  It is not required to enter the pool again.  Furthermore, any flights behind it in the queue must wait until it has landed before landing themselves, even if they have earlier unconstrained landing times.

\text{ }\indent The essence of our decision-making problem is to decide the release times $R_i$ for flights $i\in\mathcal F$.  At any time $t\in\mathcal T$, we are aware of the following dynamic information:\\
\begin{itemize}
\item The latest ETAs $X_i(t)$ of all flights $i\in\mathcal F$ that have not yet landed, including those that are in the queue, those that are in the pool and those that are yet to arrive in the pool;
\item The ordering of any planes that are in the queue;
\item Any weather forecast information that remains relevant at time $t$, which may be summarized by $\{T_0(t),\;T_1(t)\}$ if $t<U_0$ and by $T_1(t)$ if $U_0<t<U_1$;
\item If a flight $j$ is currently `in service', meaning that $t\geq A_j$ but $t<L_i+M_{ij}$ (where $i$ is the predecessor of $j$ in the queue), then we are aware of how much time $j$ has spent in service so far.\\
\end{itemize}
\text{ }\indent The above information comprises the \emph{system state} in our problem.  We are also aware of the weight classes $w_i$ and cost parameters $g_i$ for all $i\in\mathcal F$.  Actions can be taken at any time $t\in\mathcal T$ for which the pool is non-empty, implying that the set of decision epochs is uncountable; in other words, the decision-maker is able to continuously update their policy according to the latest available information, which also evolves continuously.  The `action' chosen at decision epoch $t$ is an ordered list (i.e. a tuple) of aircraft to release from the pool.  Passive actions (empty tuples) are likely to be chosen in cases where it may be advantageous to wait before deciding which aircraft to release next.  As an example, suppose the set of flights in the pool at time $t$ is $\{2,5,6,7,9\}$, where these flights have been indexed according to their positions in the original landing schedule.  Then possible actions at time $t$ include $(6,9,2,7,5)$, $(6,7)$, $(2)$ and $\emptyset$.

\text{ }\indent It remains for us to specify our cost mechanism and objective function.  We suppose that, in making decisions about when to release aircraft from the pool (and thereby determining the runway landing sequence), we have two broad objectives in mind: (i) flights should be able to land as near as possible to their scheduled landing times; (ii) flights should have their total flight times minimized.  These two objectives are not entirely unrelated, but they are quite different in nature.  The first objective is related to punctuality of air transport operations and the need to avoid disruption to airline schedules (taking into account the transfer times needed between aircraft flight legs, etc.).  The second objective is based on the need to avoid overly long flight times, which would compromise safety and increase fuel costs.  Similar versions of both objectives have been considered in many previous studies (\cite{Bennell2017}).

\text{ }\indent Specifically, for flight $i\in\mathcal F$, we define
\begin{align}C_i^{[S]}=\left(\left[L_i-a_i-\gamma^{[S]}\right]^+\right)^2,\;\;\;\;\;\;C_i^{[W]}=\left(\left[L_i-(A_i-\rho_i)-\gamma^{[W]}\right]^+\right)^2,\label{obj_defns}\end{align}
where $\gamma^{[S]}\geq 0$ and $\gamma^{[W]}\geq 0$ are parameters given as model inputs.  The interpretation is that, in the first objective, a delay of up to $\gamma^{[S]}$ time units between the actual landing time and the scheduled landing time is `acceptable', but beyond this we incur a penalty which increases quadratically with the amount of delay.  To interpret the second objective, note that 
$A_i-\rho_i$ is the time that flight $i$ would land on the runway if it were released from the pool immediately and did not incur any queueing delays.  We calculate the amount of `extra airborne time' incurred as the difference between the actual landing time, $L_i$, and $A_i-\rho_i$.  This `extra airborne time' can be incurred as a result of being held in the pool \emph{or} being delayed in the queue, so it is not necessarily the case that keeping flight $i$ waiting in the pool leads to an increase in $C_i^{[W]}$ (indeed, in many situations, keeping a flight waiting in the pool implies that it will spend less time waiting in the queue after its release; this is discussed further in the next section).

\text{ }\indent An advantage of our simulation-based solution approach is that the cost functions $C_i^{[S]}$ and $C_i^{[W]}$ can be made almost arbitrarily complicated without affecting tractability; they do not necessarily need to take quadratic or even polynomial forms, for example.  By considering only the positive part of $L_i-a_i-\gamma^{[S]}$ we avoid penalizing early landing times, but penalties for earliness could easily be incorporated into our model if needed.  We also note that there is an obvious choice available for $\gamma^{[S]}$, as 15 minutes is often regarded as a threshold for the purpose of classifying delays in the aviation industry (\cite{Ball2010,Belcastro2018}).

\text{ }\indent The objective of the problem is to minimize the dual-criteria objective function
\begin{equation}\sum_{i\in\mathcal F}g_i\times \left[\theta^{[S]} C_i^{[S]}+\theta^{[W]} C_i^{[W]}\right],\label{objective}\end{equation}
where $g_i\in[0,1]$ is the relative cost parameter for flight $i$ and $\theta^{[S]}$ and $\theta^{[W]}$ are positive-valued weights, normalized so that $\theta^{[S]}+\theta^{[W]}=1$.  The value of $g_i$ might depend on the number of passengers carried or the potential knock-on effects to other flights if the landing is delayed, for example. 

\text{ }\indent Traditionally one would approach a stochastic, dynamic optimization problem by aiming to find an optimal `policy', mapping states to actions.  In this case we have a problem with a high-dimensional, continuous state space which is obviously beyond the scope of exact solution by approaches such as dynamic programming.  We propose to take a simheuristic approach and aim to use simulation methods to update our `belief' of the optimal landing sequence as time (and uncertainty) evolves.  Further details are provided in the next section.\\

\section{Solution methodology}\label{sec:methodology}

\text{ }\indent The main solution approach of interest in our study is a simheuristic approach\footnote{See \cite{Juan2015} for a useful review of simheuristic methods for stochastic optimization problems.}, implemented in real time, in which we continuously update performance estimates for a small number of `candidate' solutions (sequences), aiming to discover new strongly-performing solutions and discard weaker ones as time progresses.  Within the taxonomy of metaheuristic algorithms, our approach may be compared to a Variable Neighborhood Search (VNS), as we occasionally force the algorithm to migrate to a different region of the solution space if it seems to be making no further progress in finding improved solutions in its current neighborhood.

\text{ }\indent We provide details of our simheuristic algorithm in Section \ref{sec31}.  In Section \ref{sec32} we describe an alternative, simpler approach, based on considering expected values of random variables rather than using simulation.  In Section \ref{sec33} we also suggest some additional policies, including the simple `first-come-first-served' rule, that can be used as benchmarks for our other heuristics.\\

\subsection{The simheuristic approach}\label{sec31}

\begin{figure}[htbp]
    \begin{center}
        \includegraphics[width=16.5cm]{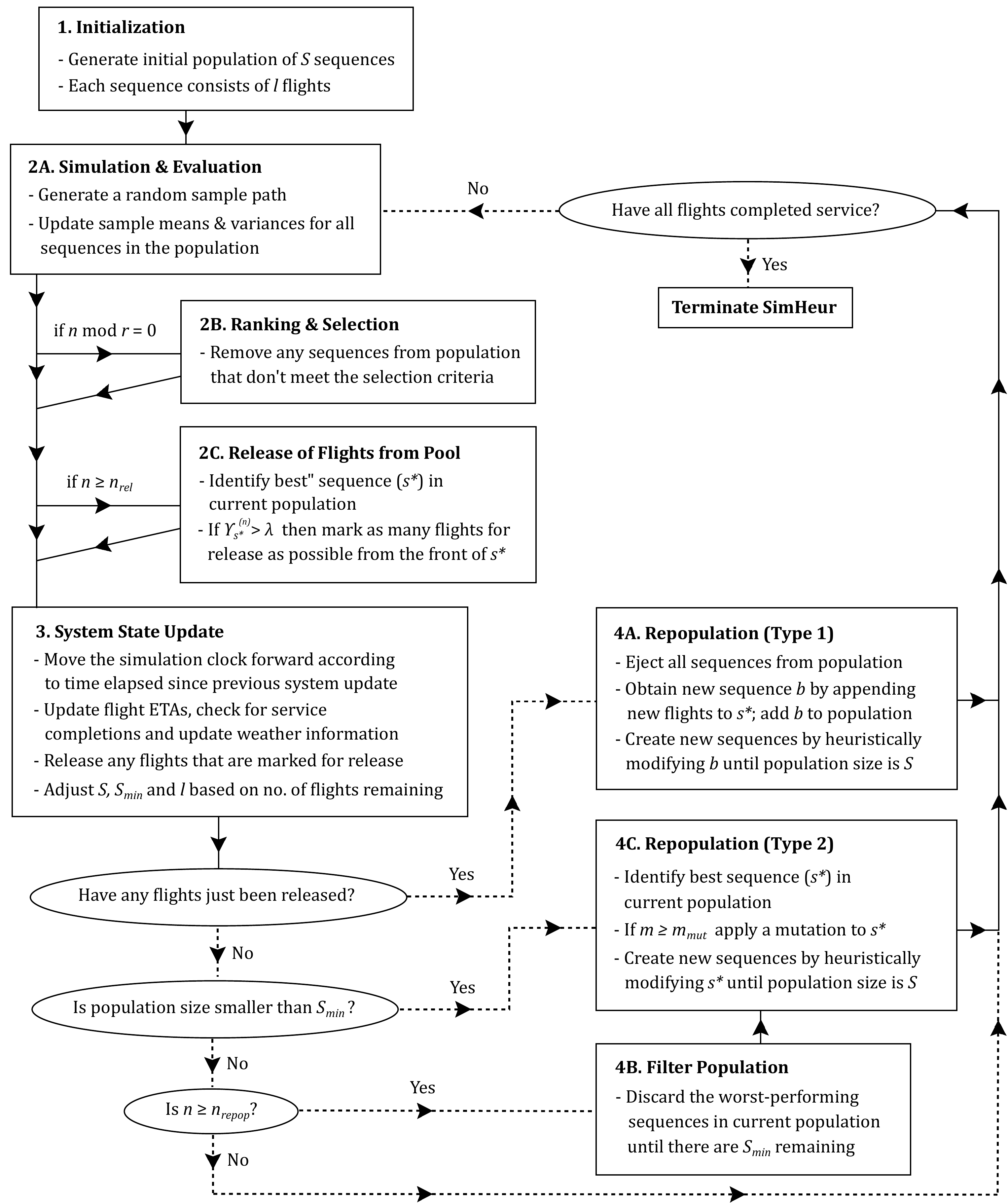}
    \end{center}
    \caption{The SimHeur algorithm.}\label{sh_fig}
\end{figure}

\text{ }\indent Let the flights in $\mathcal F$ be indexed in ascending order by their $a_i+\Delta_i^{\text{pre}}$ values, which are assumed known at the beginning of $\mathcal T$.  Hence, our initial belief is that flight $1$ will be the first to arrive in the terminal area and flight $F:=|\mathcal F|$ will be the last.  Figure \ref{sh_fig} shows an outline of the main steps in our simheuristic algorithm, which we refer to as `SimHeur' for convenience.  In the remainder of this subsection we describe these steps in more detail.\\

\textit{Step 1: Initialization}

\text{ }\indent We initialize a set (or `population') $\mathcal S_0$ of solutions (sequences), where each sequence $s\in\mathcal S_0$ is a tuple of length $l\in\mathbb{N}$ and specifies the next $l$ flights that will land on the runway, in order, if this sequence is followed.  It will be convenient to let $\mathcal S_t$ denote the population at time $t\in\mathcal T$.  The number of sequences in $\mathcal S_0$ is denoted by $S\in\mathbb{N}$ and it will also be the case that $|\mathcal S_t|\leq S$ for all $t\in\mathcal T$.  Let $S_{\min}\leq S$ be an integer that represents the minimum allowable population size.  This means that if $|\mathcal S_t|<S_{\min}$ at some time point $t$ then we must `repopulate' $\mathcal S_t$ by generating new sequences until its size is restored to $S$.  We also initialize an iteration counter, $n=0$, and an additional counter $m=0$.  The parameters $S$, $S_{\min}$ and $l$ remain constant until the final stages of SimHeur's running time, at which point it becomes necessary to reduce their values in order to ensure that it is still feasible to find $S$ distinct $l$-tuples consisting of only the flights in $\mathcal F$ that haven't already been released.  Indeed, after all flights in $\mathcal F$ have been released it is necessary to set $S=S_{\min}=l=0$, at which point steps 2A-2C and 4A-4C in Figure \ref{sh_fig} become redundant, but the algorithm continues to run until all flights in the queue have completed service.

\text{ }\indent In our numerical experiments in Section \ref{sec:numerical}, we generate the initial sequences $s\in\mathcal S_0$ by making random changes to a `first-come-first-served' sequence, as follows:\\

\begin{enumerate}[(a)]
\item Let $\mathcal S_0$ consist of only one sequence, $(1,2,...,l)$.
\item Make a change to the sequence $(1,2,...,l)$ by applying a heuristic move operator, denoted $H$.  If the new sequence is not already in $\mathcal S_0$, then add it.  Otherwise, repeat this substep.
\item If $|\mathcal S_0|=S$, terminate.  Otherwise, return to substep (b).\\
\end{enumerate}

Details of the heuristic move operator $H$ are provided in Appendix \ref{AppA}.\\

\textit{Step 2A: Simulation and Evaluation}

\text{ }\indent In this step we update performance estimates of all sequences in our current population $\mathcal S_t$.  For each flight $i$ in a particular sequence $s\in\mathcal S_t$ we must estimate its contribution to the objective function (\ref{objective}) given that $s$ is followed.  We do this by randomly sampling a sequence of events (referred to as a \emph{sample path} or \emph{sample trajectory}), denoted by $\Omega_t$.  The sample path $\Omega_t$ includes predictions, denoted by $\tilde{Q}_i$, of pool arrival times for all flights $i\in\mathcal F$ yet to arrive in the pool (i.e. $t<Q_i$), and also predictions of times needed for the remaining part of the journey for flights that are either still in the pool (i.e. $Q_i<t<R_i$) or still enroute to the runway (i.e. $R_i<t<A_i$).  It also includes predictions $\tilde{M}_{ij}$ of the separation times for all consecutive pairs $(i,j)$ of aircraft in the queue, and (if relevant) predictions $\tilde{U}_0$ and $\tilde{U}_1$ for the starting and ending times of any future period of bad weather.  We assume that all of these values are sampled from the correct distributions described in Section \ref{sec:formulation}.  For example, the predictions $\tilde{Q}_i$ can be obtained by numerically integrating the Inverse Gaussian density function 
$$f_{Q_i}(u\;|\;X_i(t))=\frac{X_i(t)-\tau-t}{\sigma_i\sqrt{2\pi u^3}}\exp\left(-\frac{(X_i(t)-\tau-t)^2}{2\sigma_i^2 u}\right)\;\;\;\;(u>0)$$
in order to obtain the distribution function of $Q_i$, and then choosing a quantile from this according to a uniform $[0,1)$ distribution (the `inverse transform sampling' method).  A similar method can also be used for sampling the remaining travel times for flights that have already arrived in the pool or been released from the pool; for example if flight $i$ is still in the pool, its remaining travel time is sampled from an Inverse Gaussian distribution with mean $\tau$, variance $\sigma_i^2\tau$.  If a service is in progress, then (as described in Section \ref{sec24}) the elapsed service time is known and we sample the remaining service time from a conditional gamma distribution.

\text{ }\indent After $\Omega_t$ has been generated, the future release times, landing times etc. of the flights in a particular sequence $s\in\mathcal S_t$ can be worked out in a deterministic way based on the generated timings of events, and this yields an overall cost estimate for sequence $s$.  For clarity, we assume here (only for the purposes of estimating costs under different possible sequences) that if a particular sequence $s\in\mathcal S_t$ is followed starting from time $t$, then the following procedure is used:\\
\begin{enumerate}[(a)]
\item If the first flight in $s$ is already in the pool, then it is released immediately.  Otherwise, we wait until it arrives in the pool, without releasing any other flights in the meantime, and then release it immediately.
\item Substep (a) is repeated for the next sequential flight in $s$, and this process is repeated until all flights in $s$ have been released.\\
\end{enumerate}

\text{ }\indent It should be noted that the cost estimate for $s$ obtained by the above procedure is based only on the $l$ flights that are included in sequence $s$; it does not take into account the pool holding times and landing times of any subsequent flights, even though these would be affected by the sequencing decisions for flights in $s$.  In this respect, our cost estimates are only based on `looking ahead' by a limited amount of time into the future, and it is important to avoid making $l$ too small in order to avoid being too myopic.  On the other hand, larger values of $l$ are associated with too much computational expense and compromise the performance of SimHeur.  In practice, the values of $l$ that we use are sufficiently small to ensure that the sequences $s\in\mathcal S_t$ consist only of flights that are already in the `tactical uncertainty' stage (i.e. $t>h_i$ for flights $i\in s$) and, hence, we do not consider sequencing options for flights that are yet to realize their pre-tactical uncertainty.  This explains why there is no loss of generality in assuming that $\Delta_i^{\text{pre}}$ is realized for all $i\in\mathcal F$ at the beginning of $\mathcal T$ and scaling the time units accordingly.

\text{ }\indent During the running of SimHeur we revisit step 2A many times and acquire cost estimates for each population member $s\in\mathcal S_t$ at many different time points (and using many different sample paths).  Let $(t_1,t_2,...)$ denote the sequence of time epochs at which these estimates are obtained, with $(\Omega_{t_1},\Omega_{t_2},..)$ being the corresponding sequence of sample paths generated.  Let $J_s^{(n)}$ denote the $n^{\text{th}}$ cost estimate for sequence $s$, given by sample path $\Omega_{t_n}$.  We note here that the iteration counter $n$ is reset to zero in some later steps of the algorithm when new sequences are added (see steps 4A-4C) and our notation assumes that any sequence included in the population at epoch $t_n$ is also included at epochs $t_j$ for $j<n$, i.e. $\mathcal S_{t_n}\subseteq \mathcal S_{t_{n-1}}\subseteq ... \subseteq \mathcal S_{t_1}$.  At time $t_n$, we update two overall performance indicators of sequence $s\in\mathcal S_{t_n}$, denoted $V_s^{(n)}$ and $W_s^{(n)}$, as follows:
\begin{align}&V_s^{(n)}:=(1-\psi_n)V_s^{(n-1)}+\psi_n J_s^{(n)},\label{rl_update1}\\
&W_s^{(n)}:=(1-\psi_n)W_s^{(n-1)}+\psi_n \left(J_s^{(n)}\right)^2.\label{rl_update2}\end{align}
\text{ }\indent (We assume $V_s^{(0)}=W_s^{(0)}=0$ for all sequences $s$.)  If we set $\psi_n:=1/n$ then equations (\ref{rl_update1})-(\ref{rl_update2}) are equivalent to simple averaging over the $n$ cost estimates (and their squares) obtained so far.  However, it might make more sense to let the sequence $(\psi_n)_{n=1}^\infty$ tend to zero at a slower rate than this, or even to set it to a small constant value, in order to obtain a reinforcement learning-style rule.  The reason for this is that each cost estimate $J_s^{(n)}$ is obtained using the latest up-to-date information (from the most recent system state update), so the more recent estimates should be more accurate and should arguably carry more weight.

\text{ }\indent There is one further performance measure that we update in this step.  Consider the flights that are already in the queue at time $t_n$ and suppose flight $i\in\mathcal F$ is the last flight in the queue, i.e. the most recent flight to have been released from the pool.  Also, let $\tilde{L}_i^{(n)}$ denote the actual landing time for flight $i$ under the sample path $\Omega_{t_n}$.  For each sequence $s\in\mathcal S_t$, let $j_s$ be the first sequential flight in $s$ and define the binary variable $\xi_s^{(n)}$ as follows:
$$\xi_s^{(n)}:=\begin{cases}1,&\text{ if }\tilde{A}_{j_s}^{(n)}> \tilde{L}_i^{(n)}+\tilde{M}_{i,j_s}^{(n)},\\
0,&\text{ otherwise,}\end{cases}$$
where $\tilde{A}_{j_s}^{(n)}$ and $\tilde{M}_{i,j_s}^{(n)}$ are (respectively) the unconstrained landing time for flight $j_s$ and the landing time separation between $i$ and $j_s$ under sample path $\Omega_{t_n}$.  Hence, if $\xi_s^{(n)}=1$, this indicates that flight $j_s$ does not arrive early enough to be able to land at the earliest possible moment after flight $i$'s landing.  One might say that there is some `idle runway time' caused by the late arrival of flight $j_s$.  We then define $\Upsilon_s^{(n)}$ as a (possibly weighted) average of $\xi_s^{(n)}$ over all sample paths, as follows:
$$\Upsilon_s^{(n)}:=\begin{cases}0,&\text{ if }n=0,\\
(1-\psi_n)\Upsilon_s^{(n-1)}+\psi_n \xi_s^{(n)},&\text{ otherwise.}\end{cases}\label{rl_update3}$$
\text{ }\indent If $\Upsilon_s^{(n)}$ happens to be very small, then this suggests that there is little benefit in releasing flight $j_s$ from the pool immediately, as it is likely to be forced to wait in the queue and will have to wait until time $L_i+M_{i,j_s}$ before it is able to land.  It may be advantageous to delay its release so that we can acquire more information (from system state updates) before making a final decision about which flight to release next from the pool.  We elaborate on this further in step 2C.\\

\textit{Step 2B: Ranking and Selection}

\text{ }\indent After $n$ cost evaluations have been performed we can rank the sequences in $\mathcal S_{t_n}$ according to their $V_s^{(n)}$ values and remove any sequences that appear to perform poorly.  In this step we utilize a standard approach in simulation optimization known as `ranking and selection' (see, for example, \cite{Nelson2013}).  First, we note that this step is not actually performed at each iteration $n\in\mathbb{N}$ in our algorithm; instead, as shown in Figure \ref{sh_fig}, it is only performed when $n$ is a multiple of $r$, for some pre-determined $r\in\mathbb{N}$.  The reason for this is that the ranking and selection process involves pairwise comparisons between all sequences in our current population, and this can be computationally expensive, so there is little value in performing this step at every iteration given that the differences in $V_s^{(n)}$ and $W_s^{(n)}$ values on consecutive iterations are likely to be small.

\text{ }\indent In this step we consider each sequence $s\in\mathcal S_{t_n}$ in the current population and retain it in the population if and only if it satisfies
\begin{equation}V_s^{(n)}\leq V_{s'}^{(n)}+Z_{s,s'}\;\;\;\;\;\;\forall s'\in\mathcal S_{t_n}\setminus\{s\},\label{rs_criterion}\end{equation}
where the threshold $Z_{s,s'}$ is given by
$$Z_{s,s'}=\left(\frac{z_{1-\eta/2}^2}{n-1} \left(W_s^{(n)}+W_{s'}^{(n)}-{[V_s^{(n)}]}^2-{[V_{s'}^{(n)}]}^2\;\right)\right)^{1/2}$$
and $z_{1-\eta/2}$ is the $(1-\eta/2)^{\text{th}}$ quantile of the standard normal distribution.  (We assume that $r$ is sufficiently large to justify using the normal distribution rather than the Student's $t$-distribution.)  We note that, since $Z_{s,s'}$ is positive, the sequence with the smallest sample mean is guaranteed to be retained in the population.\\

\textit{Step 2C: Release of Flights from Pool}

\text{ }\indent Let $s^*$ denote the sequence in our current population $\mathcal S_t$ with the smallest value of $V_s^{(n)}$ after $n$ cost evaluations have been performed.  If $n\geq n_{\text{rel}}$, where $n_{\text{rel}}\in\mathbb{N}$ is a pre-determined threshold, then we check to see whether $\Upsilon_{s^*}^{(n)}$ exceeds another pre-determined value $\lambda\in[0,1)$.  If the additional condition $\Upsilon_{s^*}^{(n)}>\lambda$ holds, then it is decided that the flights at the front of sequence $s^*$ should be released as soon as possible if they are already in the pool.

\text{ }\indent Specifically, suppose the conditions $n\geq n_{\text{rel}}$ and $\Upsilon_{s^*}^{(n)}>\lambda$ hold and let
$$u:=\max\{j\in\mathbb{N}\;:\;Q_i\leq t\text{ for all flights }i\text{ in the first }j\text{ positions of }s^*\}.$$
\text{ }\indent In other words, $u$ is the number of positions in sequence $s^*$ that we are able to count, starting from the beginning, without getting to an aircraft that isn't in the pool yet.  If $u\geq 1$, then we should release these $u$ aircraft as soon as possible, so that they join the queue in the same order that they appear in $s^*$.  As shown in Figure \ref{sh_fig}, these $u$ aircraft are only `marked' for release at this stage.  They are actually released in step 3, following the next system state update, in order to ensure that the state information (including the weather state, for example) at their time of release is accurate.  On the other hand, if $u=0$, then no aircraft should be released.

\text{ }\indent We note that, as with several other parameters in our algorithm, setting the value of $n_{\text{rel}}$ involves a `trade-off'.  Larger values enable us to be more confident that the sequence $s^*$ is genuinely the best sequence available due to the greater number of cost evaluations performed, but by requiring $n$ to be large before any flights are released, we might delay their release for too long and incur greater costs as a result.  Furthermore, as explained in step 2A, the condition $\Upsilon_{s^*}^{(n)}>\lambda$ is designed to ensure that we can derive some benefit from delaying the release of a flight from the pool in situations where the flight is likely to be delayed in the queue anyway (and therefore an early release would not imply an earlier landing time).  The benefit of delaying the release is that we are able to acquire more information (through system state updates) before deciding which flight should be committed to the queue next.  However, our numerical experiments indicate that $\lambda$ should be set to a very small value in order to give the best results, and indeed $\lambda=0$ is often the best choice.  The reason for this is that if any unnecessary `idle runway time' occurs (i.e. a flight arrives at the runway later than its earliest feasible landing time based on time separations), this can cause delays to many subsequent flights, implying a very significant increase in the value of the objective function (\ref{objective}).  Therefore, if there is even the slightest possibility of idle runway time occurring, we may wish to release the next flight as soon as possible.  Setting $\lambda=0$ effectively implies that we only need to find one random sample path (among possibly thousands) with unnecessary idle runway time in order for the condition $\Upsilon_{s^*}^{(n)}>\lambda$ to be met.\\

\textit{Step 3: System State Update}

\text{ }\indent All of the steps in our algorithm require some computational effort.  In this step we move the simulation `clock' forward according to the amount of time elapsed since the previous system state update (or since initialization, if this step is being encountered for the first time) and update the latest system state information.  Specifically, if $\delta t$ is the amount of time elapsed since the previous update, then the current time $t$ should be incremented by an amount proportional to $\delta t$.  It is then necessary to update the ETAs $X_i(t)$ for all flights $i\in\mathcal F$ such that $t<A_i$.  In addition, we need to check whether any service phase completions have occurred during this time increment and (if necessary) also update the weather forecast and the current weather state.

\text{ }\indent In our computer implementation, we have found that the most efficient approach is to pre-generate all random events and their timings, so that we have a pre-generated `actual' sample path $\Omega^{\text{true}}$ (hidden from the decision-maker).  
Then, in order to update the system state at a new time point, we only need to `look up' the relevant information in $\Omega^{\text{true}}$ rather than sampling from any distributions again.  Further details about how this is done can be found in Appendix \ref{AppB}.

\text{ }\indent This step also involves releasing any flights that have been marked for release (see step 2C), so that new flights are added to the queue if necessary.  Also, as noted in step 1, we may need to reduce the values of $S$, $S_{\min}$ and $l$ in this step if there are only a few flights remaining that haven't already been released from the pool.\\

\textit{Step 4A: Repopulation (Type 1)}

As shown in Figure \ref{sh_fig}, this step follows step 3 in the case where at least one flight has just been released from the pool.  In this case, the newly-released flights are no longer eligible to be included in the sequences in our population (as we only consider sequencing decisions for flights that haven't been released yet), so we perform a `reset' by removing all sequences from our current population, erasing the information $\{V_s^{(n)},W_s^{(n)},\Upsilon_s^{(n)}\}$ for all removed sequences and setting $n=0$.  We then create a new population $\mathcal S_t$ by performing the following substeps:\\
\begin{enumerate}[(a)]
\item Form a new sequence $b$ consisting of the flights in positions $u+1,u+2,...,|s^*|$ of the sequence $s^*$ that was chosen as the best sequence in the previous population.  (Note that $u$ represents the number of flights that have just been released, as defined in step 2C.)
\item Consider all flights $i$ that have not yet been added to the queue and are not already included in $b$.  Among these flights, select the one with the earliest ETA $X_i(t)$ and append it to the end of $b$.
\item Repeat substep (b) as many times as necessary until the number of flights in $b$ is $l$.
\item Generate an extra $S-1$ sequences, where each new sequence is formed by applying the heuristic move operator $H$ described in Appendix \ref{AppA} to the sequence $b$ formed above, in order to obtain a new population of $S$ sequences.\\
\end{enumerate}

\text{ }\indent We note that an alternative approach in this step would be to retain all of the sequences that were in the previous population, with some extra flights appended in order to increase their lengths to $l$, and also retain the information $\{V_s^{(n)},W_s^{(n)},\Upsilon_s^{(n)}\}$ for these sequences in order to avoid `throwing away' potentially useful information.  However, after removing flights from a particular sequence and adding new ones, the subsequent cost estimates for that sequence are likely to be very different from those obtained before flights were removed/added, so the `old information' is likely to be of limited use and may actually hinder the acquisition of accurate cost estimates for the amended sequence.  In our experiments, we have found that performing a `complete population reset' as described in substeps (a)-(d) above tends to be a more efficient approach.

\text{ }\indent After the new population $\mathcal S_t$ has been created, the algorithm returns to step 2A (simulation and evaluation).\\

\textit{Step 4B: Filter Population}

\text{ }\indent Following the latest system state update (step 3), it may be beneficial to add new sequences to the population $\mathcal S_t$ and evaluate these (in step 2A) according to the latest state information.  In order to do this, we need to `make room' for the new sequences by removing some of the weaker sequences from $\mathcal S_t$.  This step is only performed if $n\geq n_{\text{repop}}$, where $n_{\text{repop}}\in\mathbb{N}$ is a pre-determined threshold, because we must have performed a sufficient number of iterations to be able to reliably judge which sequences in $\mathcal S_t$ are the weakest.

\text{ }\indent In this step we simply rank the sequences in $\mathcal S_t$ according to their sample means $V_s^{(n)}$ (comparable to the `fitness' estimates used in metaheuristic algorithms) and remove the $(S-S_{\min})$ sequences with the highest values, so that the new population size is $S_{\min}$.  If $S_{\min}=1$ then we only retain the `best' sequence in the current population, but there are some potential advantages in choosing a larger $S_{\min}$ value: (i) the sample means $V_s^{(n)}$ are associated with sampling error and the ranking order of the sequences may change as we repeat step 2A more times; (ii) even if there was no sampling error, the system state is continuously evolving and new information might imply that the ranking order should be changed.  Note that sequences can be removed in this step even if they satisfy the `ranking and selection' criterion (\ref{rs_criterion}) given in step 2B.

\text{ }\indent After this step, we then add new sequences to the population in step 4C.\\

\textit{Step 4C: Repopulation (Type 2)}

\text{ }\indent As shown in Figure \ref{sh_fig}, this step can be reached in two different ways.  If $n\geq n_{\text{repop}}$, then we reduce the population size to $S_{\min}$ as described in step 4B before arriving at this step.  On the other hand, we might also reach this step if the population size has been reduced to $S_{\min}$ (or smaller) following the `ranking and selection' process in step 2B.  In either case, we can assume that the current population size $|\mathcal S_t|$ is not greater than $S_{\min}$.

\text{ }\indent The purpose of this step is to add new sequences to the population until its size is restored to $S$.  Recall that $m$ is a counter, initialized with a value of zero in step 1.  We assume that $m_{\text{mut}}\in\mathbb{N}$ is a pre-determined threshold which determines the point at which we should `mutate' the best sequence in our current population.  The substeps are described below.  \\

\begin{enumerate}[(a)]
\item Let $s^*$ denote the sequence in the current population $\mathcal S_t$ with the smallest sample mean $V_s^{(n)}$.
\item If $s^*$ was already included in the population when we last entered step 4C, then increase $m$ by 1.  Otherwise, set $m=0$.
\item If $m\geq m_{\text{mut}}$, create a new sequence $b$ by applying a random mutation to $s^*$ and then set $m=0$.  Otherwise, set $b=s^*$.
\item Make a change to the sequence $b$ by applying a heuristic move operator, denoted $H$.  If the new sequence is not already in $\mathcal S_t$, then add it.  Otherwise, repeat this step.
\item If $|\mathcal S_t|=S$, terminate.  Otherwise, return to step (d).\\
\end{enumerate}

\text{ }\indent The heuristic move operator $H$ is the same one used in steps 1 and 4A, and is described in Appendix \ref{AppA}.  The `random mutation' referred to in substep (c) above is described in Appendix \ref{AppC}.  We note that the purpose of the mutation is to force the algorithm to explore a different part of the solution space, and is consistent with the VNS methodology.  The condition $m\geq m_{\text{mut}}$ indicates that the algorithm has been through the repopulation process $m$ times without successfully finding a new sequence that performs better than $s^*$.  This suggests that a local optimum has been found, and hence we should migrate to another region of the solution space.

\text{ }\indent We also reset $n$ to zero in this step and erase the information $\{V_s^{(n)},W_s^{(n)},\Upsilon_s^{(n)}\}$ for sequences that have been retained from our previous population.  The reason for this is that we wish to compare the retained sequences with newly-added sequences, and in order to ensure a fair comparison it is best to ensure that the cost estimates for all sequences (old and new) are obtained using the same set of random sample paths.  The algorithm then returns to step 2A (simulation and evaluation).\\

\text{ }\indent We have now described all steps in the SimHeur algorithm.  As shown in Figure \ref{sh_fig}, the algorithm terminates when all flights in $\mathcal F$ have completed service, at which point we obtain a value for the objective function (\ref{objective}).  In the next subsection we describe `DetHeur', which operates in a more simple way by evaluating the performances of candidate sequences using expected values rather than simulated sample paths.\\

\subsection{An alternative approach based on expected values}\label{sec32}

\text{ }\indent In order to evaluate the performance of the SimHeur algorithm described in Section \ref{sec31} we will compare it to some alternative approaches.  One such alternative is to estimate the costs associated with different sequences by assuming that random variables such as $Q_i$, $M_{ij}$ etc. conform to their expected values, conditioned on the latest system state information available.  Computationally, this is much less demanding than generating random sample paths and updating the summary information $\{V_s^{(n)},W_s^{(n)},\Upsilon_s^{(n)}\}$ as described in step 2A of SimHeur, but it is also less accurate.  We refer to this simpler algorithm as `DetHeur' and the steps involved are similar to those described in Section \ref{sec31} for SimHeur, except for the following changes:\\
\begin{itemize}
\item In step 2A, we generate the sample path $\mathcal S_t$ by setting $\tilde{Q}_i=X_i(t)-\tau$ for all flights $i$ that haven't arrived in the pool yet (i.e. $t<Q_i$).  After $i$ has been released (i.e. for $t\geq R_i$) we use $\tilde{A}_i=X_i(t)$ as a deterministic prediction of its runway arrival time.  We also set $\tilde{M}_{ij}=e_{ij}$ for all consecutive pairs of flights $(i,j)$ in a particular sequence and set $\tilde{U}_0=T_0(t)$ and $\tilde{U}_1=T_1(t)$ if there is a future period of bad weather expected.  The sample mean $V_s^{(n)}$ is replaced by the single cost evaluation $J_s^{(n)}$, so that any previous cost evaluations for sequence $s$ are discarded.  We do not require $W_s^{(n)}$ or $\Upsilon_s^{(n)}$.
\item Step 2B (ranking and selection) is omitted.
\item In step 2C, the condition $n\geq n_{\text{rel}}$ is no longer required (equivalently, we might say that $n_{\text{rel}}=1$).  The condition $\Upsilon_s^{(n)}>\lambda$ is also no longer applicable.
\item The condition $n\geq n_{\text{repop}}$ is no longer required in order to enter step 4B.  Equivalently, we might say that $n_{\text{repop}}=1$.\\
\end{itemize}

\text{ }\indent All other steps are the same as in SimHeur.  Like SimHeur, the DetHeur algorithm always uses the latest available system state information (obtained in step 3) in order to make decisions, but its decisions are made based on deterministic forecasts.  An important advantage of DetHeur is that it enters the repopulation step 4C much more frequently than SimHeur (due to the removal of the condition $n\geq n_{\text{repop}}$) and this implies that it is able to spend more time searching the solution space than SimHeur, with more random mutations.  However, we conjecture that this advantage diminishes as the amount of computational power increases, because the SimHeur algorithm only needs to be able to explore `enough' of the solution space to be able to find the best solution.  If the amount of computational effort spent on exploration of the solution space is already sufficient to find the best solution, then there is no advantage to be gained by exploring further.

\text{ }\indent We note that the removal of the condition $\Upsilon_s^{(n)}>\lambda$ implies that DetHeur takes a somewhat conservative approach by releasing the first flight in $s^*$ as soon as it arrives in the pool, rather than delaying its release in order to acquire more information.  Indeed, since DetHeur does not generate random sample paths, it cannot estimate the probability of idle runway time in the same way as SimHeur.  One might argue, however, that it should release the first flight (say flight $j$) from $s^*$ if and only if the runway arrival time $\tilde{A}_j=t+\tau$ under the expected value trajectory at time $t$ satisfies $\tilde{A}_j\geq \tilde{L}_i+\tilde{M}_{ij}$, where $i$ is the latest flight to be released; in other words, the flight should not be released if we expect it to be delayed in the queue.  This seems a reasonable suggestion, but it is easy to show using experiments that DetHeur performs extremely poorly under such a rule.  As noted in Section \ref{sec31}, any unnecessary idle runway time tends to increase the objective function value very significantly due to the `knock-on' effect of one delayed landing causing another.  If we release flight $j$ at the point where $\tilde{A}_j=\tilde{L}_i+\tilde{M}_{ij}$ then there is roughly a 50\% chance that unnecessary idle runway time will occur, and this must be avoided.  We therefore opt for the conservative approach of releasing flight $j$ as soon as possible under DetHeur.

\text{ }\indent We note that in the stochastic programming literature, the idea of comparing the optimal solution to a stochastic program with the optimal solution to a corresponding `expected value' problem is widely used (see \cite{Birge2011}).  One can compare the performances of the two solutions under stochastic conditions in order to evaluate the `value of the stochastic solution' (VSS).  A similar approach can be used in our problem, by comparing the performances of SimHeur and DetHeur in order to investigate the benefits of being able to simulate random events based on knowledge of the underlying probability distributions.\\

\subsection{Other sequencing strategies and benchmarks}\label{sec33}

\text{ }\indent In our numerical experiments in Section \ref{sec:numerical} we make use of some additional benchmarks for evaluating the performances of SimHeur and DetHeur.  The first of these is the cost associated with a simple `first-come-first-served' (FCFS) rule, in which flights are released immediately when they arrive at the pool ($R_i=Q_i$ for $i\in\mathcal F$).  The FCFS rule avoids pool-holding delays, but it pays no attention to scheduled landing times or the time separations required between different weight class combinations, so the resulting landing sequence may be far from optimal.

\text{ }\indent We also consider another policy obtained from a static, deterministic optimization procedure, referred to as `DStat' for short.  The DStat policy can be expressed in the form of a single sequence, i.e. an ordering of the flights $i\in\mathcal F$, which is computed \emph{once} at the beginning of the interval $\mathcal T$ (assuming knowledge of the pre-tactical delays $\Delta_i^{\text{pre}}$) and not updated at any future time points.  More specifically, at the beginning of $\mathcal T$, we assume that all random variables, including unconstrained arrival times, service times and weather transitions conform to their expected values and then aim to find the runway sequence that optimizes the objective function (\ref{objective}) under such conditions.  Although the resulting optimization problem is static and deterministic, it still has very high combinatorial complexity, and we therefore aim to solve the problem heuristically by making successive improvements to a first-come-first-served sequence until no further improvements can be found; further details are given in Appendix \ref{AppD}.  After the DStat sequence has been obtained, its performance under the stochastic conditions of our model can be evaluated using the `actual' values of the random variables stored within the sample path $\Omega^{\text{true}}$.

\text{ }\indent Although the DStat policy is similar to DetHeur in that both algorithms treat the problem as deterministic and rely upon expected values, the DStat policy is much more simplistic than DetHeur as it is not a \emph{dynamic} policy; that is, it lacks the ability to update sequencing decisions during $\mathcal T$ in response to the latest observed information.\\

\section{Data acquisition and model calibration}\label{sec:data}

\text{ }\indent In this section we describe our use of on-time performance data for arrivals at London Heathrow Airport to estimate parameter values for the model described in Section \ref{sec:formulation} and design an appropriate flight schedule for model testing purposes.

\text{ }\indent As mentioned in Section \ref{sec:formulation}, Heathrow Airport usually operates with one of its two runways used exclusively for arrivals.  We decided to select a particular day of operations and look up historical on-time performance data for the arriving flights scheduled during a particular part of that day.  Our chosen day was August 1st, 2019 and we considered the eight-hour period from 6:00AM (inclusive) to 2:00PM (non-inclusive), during which there were 323 arriving flights scheduled.  We take these 323 flights as our set of flights $\mathcal F$.  We note that our chosen date ensures that the flight schedule and historical dataset are not affected by the disruption caused by the Covid-19 pandemic.  

\begin{figure}[htbp]
    \begin{center}
        \includegraphics[width=15.5cm]{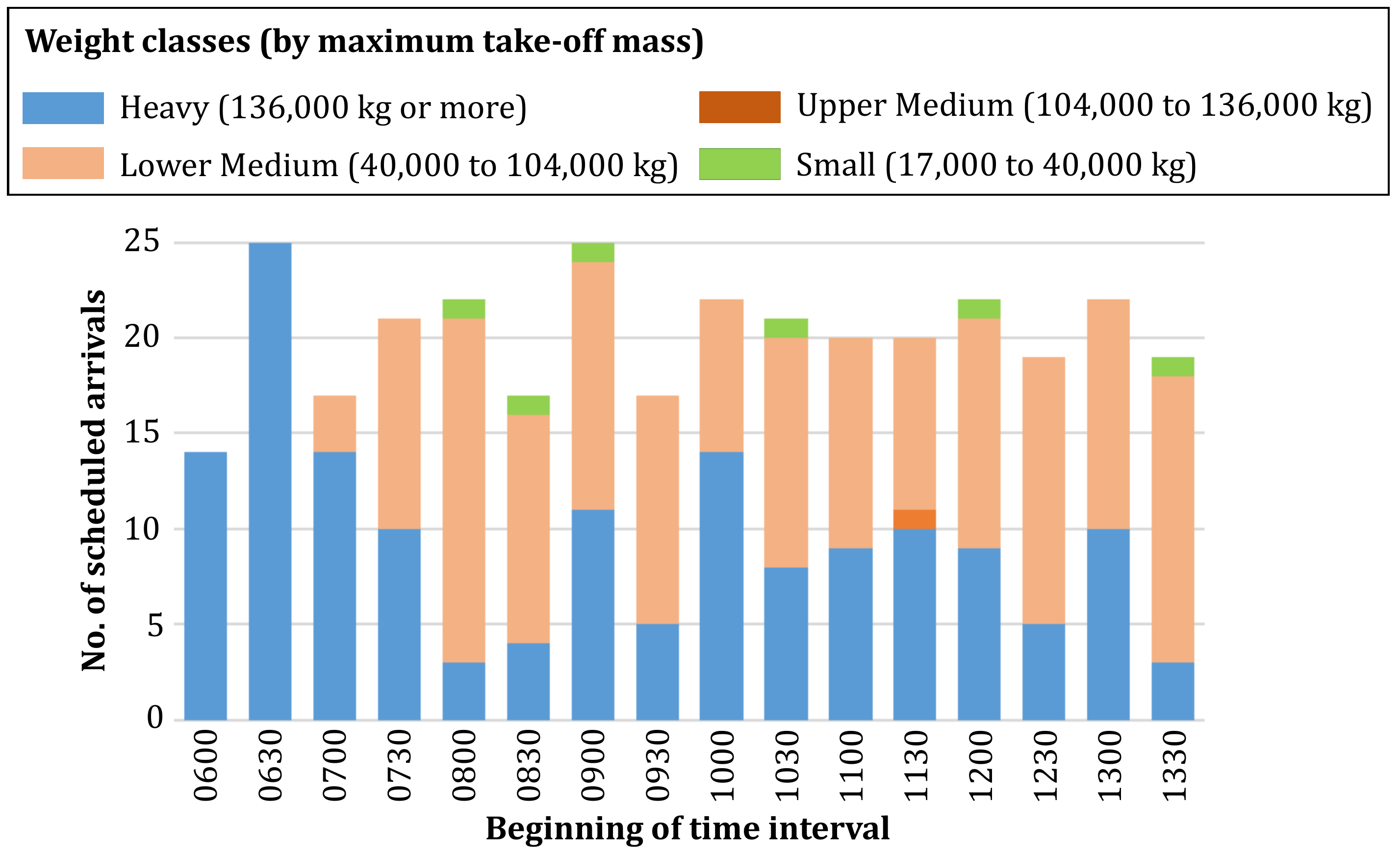}
    \end{center}
    \caption{Numbers of arrivals scheduled by half-hour time interval between 6:00AM and 2:00PM at Heathrow Airport on August 1st, 2019.}\label{sec4_fig}
\end{figure}

\text{ }\indent Figure \ref{sec4_fig} shows the numbers of arrivals scheduled within each half-hour interval during the period 6:00AM-2:00PM, with a breakdown of aircraft weight classes in each interval also provided.  Aircraft can be divided into different weight classes according to their `maximum take-off mass' (MTOM), measured in kilograms (\cite{AIS2019}).  We found that the vast majority of aircraft arriving at Heathrow on the day of interest belonged to the `heavy' and `lower medium' categories, with only a small minority belonging to the `upper medium' and `small' categories.  Larger aircraft tend to be used for long-haul, intercontinental flights, whereas smaller aircraft usually arrive at Heathrow from other European airports.  This explains why the proportion of `heavy' aircraft tends to be greatest in the earliest intervals (6:00-8:00).  Short-haul European flights are unlikely to be scheduled for arrival during these early intervals, as they would need to depart during the night in order to meet such early scheduled arrival times.

\text{ }\indent In order to specify the required time separation $e_{ij}$ between a leading aircraft of type $i$ and a following aircraft of type $j$ (where $i,j\in\mathcal W=\{\text{Heavy, Upper Medium, Lower Medium, Small}\}$) we use the matrix of separation times from \cite{Bennell2017}, which has been calibrated according to observed data from Heathrow Airport.  This matrix is shown in Table \ref{sep_times}.

\begin{table}[htbp]
\caption{Required separation times in seconds for different leader-follower aircraft weight class pairs (H = Heavy, UM = Upper Medium, LM = 
Lower Medium, S = Small).}
\centering
\begin{tabular}{cc|cccc|}
\cline{3-6}
                            &    & \multicolumn{4}{c|}{Follower} \\
                            &    & H     & UM    & LM    & S     \\ \hline
\multicolumn{1}{|c}{}       & H  & 97    & 121   & 121   & 145   \\
\multicolumn{1}{|c}{Leader} & UM & 72    & 72    & 97    & 97    \\
\multicolumn{1}{|c}{}       & LM & 72    & 72    & 72    & 72    \\
\multicolumn{1}{|c}{}       & S  & 72    & 72    & 72    & 72    \\ \hline
\end{tabular}
\label{sep_times}\end{table}

\text{ }\indent For each of the 323 arrivals at Heathrow during the eight-hour period of interest, we used historical data available at \url{http://www.flightradar24.com} to find the exact arrival times (to the nearest minute) of all flights with the same flight number (implying the same flight carrier, origin airport and destination airport) over the 360-day period from August 8th, 2018 to August 2nd, 2019.  Hence, for each flight in $\mathcal F$, we have a set of historical landing times collected over a 360-day period and are able to estimate means, variances etc. of the flight's punctuality with respect to its scheduled arrival time.  The majority of flights in $\mathcal F$ (221 out of 323) operated on at least 300 days during this 360-day period, and the total number of records we have (where the term `record' here refers to the recorded landing time of a flight in $\mathcal F$ during the 360-day historical period) is 98,814, equating to about 306 per flight on average.

\text{ }\indent Ideally, we would like to use these historical data to estimate probability distributions for the pre-tactical and tactical delays affecting the unconstrained landing times $A_i$ for flights $i\in\mathcal F$.  However, the records in our data are \emph{actual} landing times (denoted $L_i$ in our model), which may be affected by queueing delays and other congestion effects, as well as poor weather.  The records in our dataset do not provide us with a means of calculating the extent to which landing times are influenced by airport congestion and other similar effects, so we must rely on an approximate method to fit distributions for the unconstrained landing times in our model.

\text{ }\indent The method we use is as follows: for each of the 360 days in our historical period we look at the sequence of actual landings that took place during the 6AM-2PM interval and separate the flights that landed into two sets.  Set $S_1$ ($S_2$) consists of flights that landed immediately after a flight that landed earlier (later) than its scheduled landing time.  We then compare the proportions of flights landing earlier than their scheduled times in sets $S_1$ and $S_2$.  If the difference in these two proportions is statistically significant at the 5\% level, then we discard all data from the whole of that particular day when fitting the distributions for the $A_i$.  The rationale for this approach is as follows: if the queueing delays on a particular day are significant, then delays are likely to propagate, in the sense that the late arrival of one flight causes the late arrival of another.  Therefore if we restrict attention to days on which the timeliness of one landing tends to be independent of the timeliness of its immediate predecessor, we can be more confident that the recorded landing times of flights are similar to their unconstrained landing times. 

\begin{figure}[htbp]
    \begin{center}
        \includegraphics[width=16.5cm]{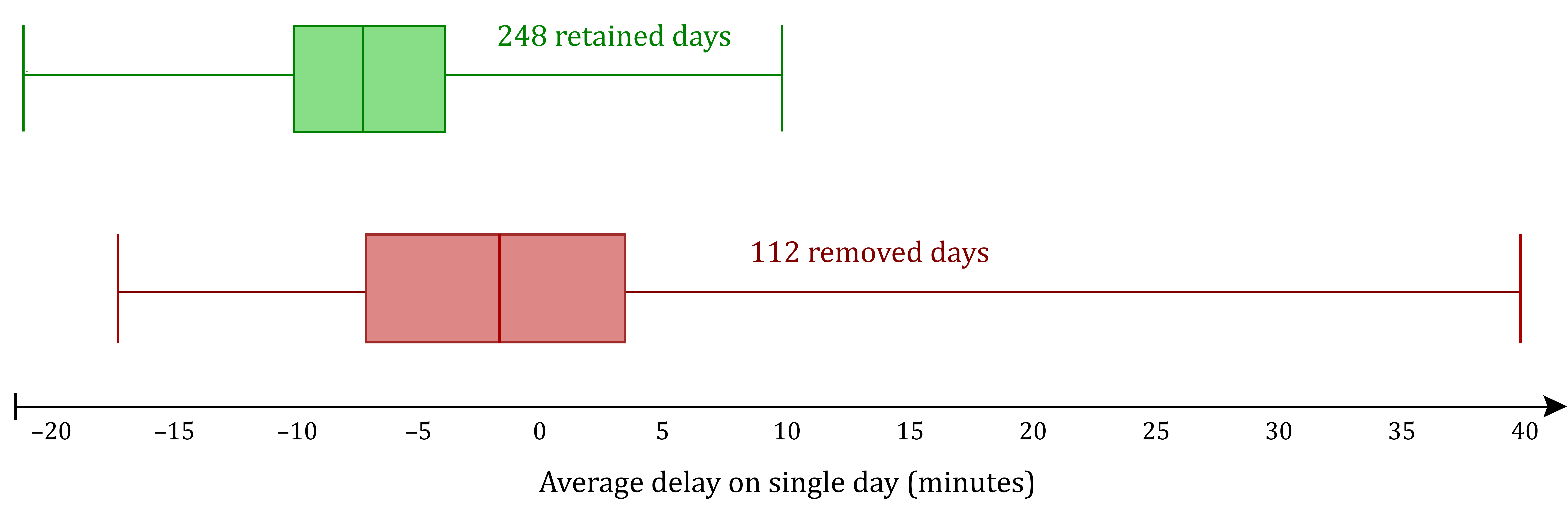}
    \end{center}
    \caption{Distributions of average delay on a single day for the 248 retained days and 112 removed days.}\label{sec4_box_plot}
\end{figure}

\text{ }\indent After carrying out the procedure described above, we identified 112 days (out of 360) on which the differences between proportions of late-arriving flights in sets $S_1$ and $S_2$ were statistically significant.  We removed all of these 112 days from our dataset and were left with 248 days' worth of data, with 199.5 records per flight on average.  The box plot in Figure \ref{sec4_box_plot} shows a comparison between the distributions of average delay over all flights in $\mathcal F$ for the 248 retained days and the 112 removed days.  (Negative delays occur when actual landing times are earlier than scheduled times.)  All of the days with average delays greater than $10.38$ minutes were among the 112 flights removed from our dataset, suggesting that our method is somewhat effective in filtering out the days with abnormally long delays that may be attributable to airport congestion.  We emphasize here that our data-filtering approach is used for estimation purposes only (i.e. to improve accuracy of estimating the $A_i$), and does not imply that we are seeking to minimize the risk of congestion occurring in our model.  Indeed, congestion occurs primarily as a result of other factors in our model, including the density of the flight schedule, required separation times, effects of bad weather etc.; it is not principally dependent on the distributions assumed for the $A_i$.

\text{ }\indent For each flight $i\in\mathcal F$, we use our 248 days' worth of data to calculate the sample mean and sample variance of the actual landing time.  We then consider how much of this variance should be explained by pre-tactical uncertainty in our model.  Recall that the mean and variance of $A_i$ are given by equations (\ref{BM_eq6a}) and (\ref{BM_eq6b}) respectively.  Our approach is to equate the expressions for $\mathbb{E}[A_i]$ and $\text{Var}(A_i)$ to the sample mean and sample variance for flight $i$ in our dataset, denoted $\bar{x}_i$ and $s_i^2$ respectively, in order to derive suitable values for our model parameters.  As mentioned in Section \ref{sec:formulation}, we rely on gamma distributions for modeling the pre-tactical delays.  Specifically, we assume that
\begin{equation}\Delta_i^{\text{pre}}=Y_i-(a_i-h_i),\label{delta_i_pre_eq}\end{equation}
where $Y_i$ is gamma-distributed with the density function
$$f_{Y_i}(t)=\frac{\beta_i^{\alpha_i}}{\Gamma(\alpha_i)}t^{\alpha_i-1}\exp(-\beta_i t)\;\;\;\;\;\;(t>0),$$
which has mean $\alpha_i\beta_i^{-1}$ and variance $\alpha_i\beta_i^{-2}$ ($\alpha_i,\beta_i>0$).  Hence, by setting $\bar{x}_i$ and $s_i^2$ equal to the expressions in (\ref{BM_eq6a}) and (\ref{BM_eq6b}) respectively, we obtain the following expressions for $\alpha_i$ and $\beta_i$:

\begin{equation}\alpha_i=\frac{(\bar{x}_i-h_i)^2}{s_i^2-\sigma_i^2(\bar{x}_i-h_i)},\;\;\;\;\beta_i=\frac{\bar{x}_i-h_i}{s_i^2-\sigma_i^2(\bar{x}_i-h_i)}.\label{alphabeta}\end{equation}
\text{ }\\
\text{ }\indent Our dataset provides values of $\bar{x}_i$ and $s_i^2$ for each $i\in\mathcal F$.  In calculating these sample statistics, we have found that it is necessary to exclude `outliers', as these tend to dominate the calculation of the sample variances $s_i^2$ and cause the distributions of $\Delta_i^{\text{pre}}$ in our model to be too platykurtic in shape.   We define an outlier as an arrival that does not occur within $120$ minutes of its scheduled arrival time.  Using this definition, about $1.06\%$ of records in our dataset are outliers.

\text{ }\indent We set the parameter $h_i$ to be 15 minutes earlier than the scheduled departure time $d_i$ in our numerical experiments, on the basis that this is a reasonable estimate for the average taxi-out time before departure (\cite{Burgain2009,Badrinath2020}).  We consider different possible cases for $\sigma_i$ in our experiments, with larger values implying a greater proportion of tactical uncertainty (as opposed to pre-tactical).  We note that the expressions in (\ref{alphabeta}) are valid only if $\sigma_i^2<s_i^2(\bar{x}_i-h_i)^{-1}$.  If we wish to consider a larger $\sigma_i^2$ value, we can simply set $\Delta_i^{\text{pre}}$ equal to its data-calibrated expected value $\bar{x}_i-a_i$ rather than sampling it using (\ref{delta_i_pre_eq}).  This represents the case where all of the variation in $A_i$ is accounted for at the tactical level.

\begin{figure}[htbp]
    \begin{center}
        \includegraphics[width=16.5cm]{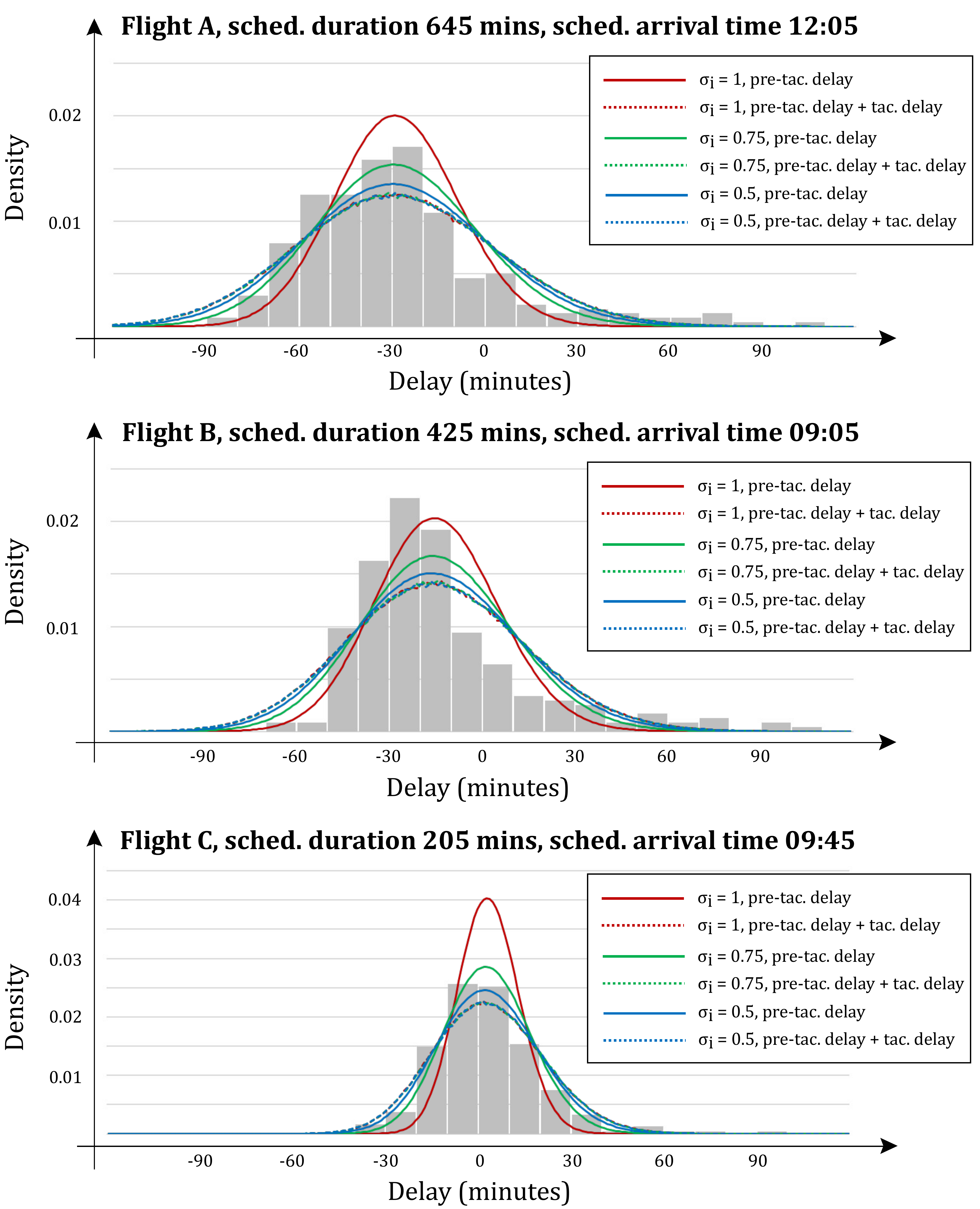}
    \end{center}
    \caption{The results of fitting distributions for $\Delta_i^{\text{pre}}$ and $(\Delta_i^{\text{pre}}+\Delta_i^{\text{tac}})$ using the empirical data, for three different flights and three different values of $\sigma_i$ in each case.}\label{sec4_data_fig}
\end{figure}

\text{ }\indent Figure \ref{sec4_data_fig} shows the results of fitting distributions for $\Delta_i^{\text{pre}}$ and $(\Delta_i^{\text{pre}}+\Delta_i^{\text{tac}})$ to the empirical data by calculating values for $\alpha_i$ and $\beta_i$ using the method described above.  For illustration purposes, we have selected 3 particular flights in $\mathcal F$ - referred to as Flights A, B and C - with noticeably different delay distributions.  The scheduled durations for A, B and C are 645, 425 and 205 minutes respectively, and the scheduled arrival times are 12:05, 09:05 and 09:45.  For each flight, the histogram shows the distribution (over all 248 days in our reduced dataset) of the delay (or `lateness') in minutes.  The red, green and blue solid curves show the distribution of the pre-tactical delay given $\sigma_i$ values of $1$, $0.75$ and $0.5$ respectively.  Similarly, the red, green and blue dashed curves show the distribution of the overall delay (including pre-tactical and tactical) for the same $\sigma_i$ values.  Recall that if $\sigma_i$ is reduced, then a larger proportion of the empirical variance is accounted for at the pre-tactical level in our model.  This explains why the solid curves become `flatter' (indicating more variance in the pre-tactical delay) as $\sigma_i$ becomes smaller.

\text{ }\indent It can also be seen that, for each of the 3 flights, the 3 dashed curves are almost indistinguishable from each other.  This shows that the value of $\sigma_i$ has almost no effect on the unconditional distribution of $A_i$.  Indeed, the role of $\sigma_i$ in our model is to determine the relative proportions of pre-tactical and tactical uncertainty affecting flight $i$'s arrival time; it does not affect the total amount of uncertainty.  However, the value of $\sigma_i$ does have a very significant effect on the nature of the decision problem formulated in Section \ref{sec:formulation}, as larger values imply that decisions must be made under higher levels of uncertainty.  This is demonstrated by the results of our numerical experiments in Section \ref{sec:numerical}.\\

\section{Numerical experiments}\label{sec:numerical}

\text{ }\indent Our computational study is based on the schedule of operations at Heathrow Airport during the eight-hour period from 6:00AM to 2:00PM on August 1st, 2019.  Our set of flights $\mathcal F$ consists of all 323 flights scheduled during this eight-hour period.  The aim of the computational study is to compare the performances of the various solution algorithms described in Section \ref{sec:methodology} and we do this by randomly generating a series of test scenarios.  In this section we begin by describing how these scenarios are generated, and then explain how time is handled in our experiments.  We then describe how the values of certain key parameters are adjusted during our study, before finally presenting a summary of results from all experiments.  

\text{ }\indent To ensure fairness, all of our experiments were performed on a single desktop computer with an Intel(R) Core(TM) i7-9700 CPU and 16GB of RAM.  Programs were implemented using Python, with the PyPy Just-in-Time compiler (\url{https://pypy.org}) used to enhance computation speed.  

\subsection{Generation of scenarios}\label{sec51}

\text{ }\indent Each test scenario consists of pre-tactical delays $\Delta_i^{\text{pre}}$ for all $i\in\mathcal F$, complete random trajectories $\{X_i(t)\}_{t\in\mathcal T}$ for all $i\in\mathcal F$, weather forecast trajectories $\{T_0(t),T_1(t)\}_{t\in\mathcal T}$ and a set of possible separation times $\{M_{ij}\}$ for each $j\in\mathcal F$ (taking into account all possibilities for the weight class of the predecessor flight $i$).  
The pre-tactical delays $\Delta_i^{\text{pre}}$ can be generated immediately by sampling from gamma distributions with parameter values configured according to the historical data as described in Section \ref{sec:data}.  The flight trajectories, weather forecasts and separation times can also be generated immediately, but this information (unlike the set of pre-tactical delays) is initially `hidden' from the decision-maker and revealed gradually as time progresses.  This hidden information comprises the `true' sample path $\Omega^{\text{true}}$ referred to in Section \ref{sec:methodology}.  Further details of how $\Omega^{\text{true}}$ is generated can be found in Appendix \ref{AppB}.

\subsection{How time is handled in our experiments}\label{sec52}

\text{ }\indent In each test scenario we wish to evaluate the performances of the SimHeur, DetHeur, FCFS and DStat algorithms, as described in Section \ref{sec:methodology}.  Recall that the SimHeur and DetHeur algorithms both operate by continuously updating their beliefs of the optimal landing sequence, with the simulation `clock' being moved forward in increments that depend on the actual amounts of time spent by these algorithms on their various computational steps (see Figure \ref{sh_fig}).  Thus, in order to allow these algorithms to perform the same number of computational steps that they would perform in a `real time' setting, we would need to run both of them for (at least) 8 hours.  However, in this computational study we wish to test a large number of different test scenarios, and it is not practical to spend 16 hours on each individual scenario.  We therefore compress the time scale in our experiments so that when the simulation clock is moved forward in step 3, the time increment added is 60 times the actual amount of CPU time elapsed since the previous increment was made.  This effectively means that one second of computation time is equated to one minute of operational time, and the SimHeur and DetHeur are only allowed to perform 1/60 of the number of steps that they would be able to perform in reality.  We conjecture that SimHeur is likely to be affected more by this time compression than DetHeur due to its higher computational demands (as discussed in Section \ref{sec32}), and hence the improvements in performance of SimHeur relative to DetHeur that we find in our experiments are likely to be `lower bounds' of the improvements that would be achievable in reality.

\text{ }\indent \textit{Note: }For clarity, when time units (e.g. minutes, hours) are referred to in the remainder of this section, these should be understood as units of `operational time' rather than `CPU time'; for example, we still refer to the time interval 6:00AM-2:00PM as being 8 hours long, although our time compression implies that it is simulated in only 8 minutes of CPU time.

\text{ }\indent We also note that this time compression does not have any effect on the performance of the FCFS algorithm, as the landing times $L_i$ for each $i\in\mathcal F$ under the FCFS policy can be calculated in a deterministic way as soon as the sample path $\Omega^{\text{true}}$ has been generated.  In practice, this means that there is no need to use a simulation clock at all when calculating the FCFS value; instead, it can be determined immediately and costs a negligible amount of CPU time.  Similarly, the DStat algorithm does not require a simulation clock and is not affected by the time compression, since its recommended sequence is obtained using only the information available at the beginning of $\mathcal T$ and its actual performance can then be calculated deterministically using the information in $\Omega^{\text{true}}$.

\text{ }\indent We count time in minutes and start counting at 5:00AM, so that $t=60$ indicates 6:00AM and $t=540$ indicates 2:00PM.  One reason for starting at 5:00AM is so that the SimHeur and DetHeur algorithms can have some initial time to evaluate possible sequences and decide on their estimated `optimal' sequences before the time window of interest (6:00AM-2:00PM) actually begins.  In addition, some early-arriving flights might arrive in the pool between 5:00AM and 6:00AM, in which case we would like to have the option of releasing them earlier than 6:00AM.  As discussed in Section \ref{sec:formulation}, the departure time $d_i$ for a particular flight $i\in\mathcal F$ is allowed to occur before time $t=0$.  In the context of our model, this implies that the time $h_i$ at which flight $i$'s trajectory $X_i(t)$ begins varying according to Brownian motion may precede $t=0$.  If this is the case, then (as discussed in Section \ref{sec24}) we must generate an initial ETA $X_i(t)$ by sampling from a Normal distribution with mean $a_i+\Delta_i^{\text{pre}}$ and variance $\sigma_i^2 |h_i|$.  Although all flights in $\mathcal F$ are scheduled to land before 2:00PM, we need to allow some extra time beyond this in order to simulate late arrivals, so we define $\mathcal T=[0,T]$ where $T=780$.  This allows for the possibility of flights landing up to 4 hours late.  (In practice, we do not need to continue simulating until time $T$ if all flights have already landed.)  Ideally, each trajectory $X_i(t)$ should vary continuously throughout $t\in\mathcal T$, but our method of pre-generating the random information (as discussed in Section \ref{sec:methodology}) implies that we can only store a finite number of $X_i(t)$ values for each $i\in\mathcal F$.  Therefore, for the purposes of these experiments, we generate and store values of each $X_i(t)$ (and also $T_0(t)$, $T_1(t)$) only for $t=0,0.01,0.02,...,T$, so that the aircraft ETAs and weather forecasts change every 0.01 minutes.

\subsection{Adjusting the values of physical parameters}\label{sec53}

\text{ }\indent In the implementation of this study, we can distinguish between `physical' parameters that affect the dynamics and costs of the problem itself and `algorithmic' parameters that are only relevant to the workings of the SimHeur and DetHeur algorithms.  The former class includes the tactical-stage variance parameters $\sigma_i$, the Erlang parameter $k$ (affecting separation time variances), the initial weather forecast parameters $t_0$ and $t_1$, the weather variance parameter $\nu$, the times $h_i$ at which flights enter the `tactical uncertainty' stage, the threshold $\tau$ that determines when sequencing decisions are made, the relative cost parameters $g_i$ for $i\in\mathcal F$, the `tolerance' parameters $\gamma^{[S]}$ and $\gamma^{[W]}$ included in the objective function (\ref{objective}) and the objective function weights $\theta^{[S]}$ and $\theta^{[W]}$.  The latter class includes parameters such as $S_{\min}$, $l$, $n_{\text{rel}}$, $n_{\text{repop}}$, $m_{\text{mut}}$ and $\lambda$ that should be chosen carefully in order for the SimHeur and DetHeur algorithms to perform well, but do not change the nature of the physical problem.  In order to avoid having too many variable factors, we have chosen to vary only the physical parameters in this computational study, and keep the values of the algorithmic parameters fixed throughout all of our experiments.  The values of the algorithmic parameters have been chosen according to a preliminary study in which the SimHeur and DetHeur algorithms were tested on some small `test' problems; further information can be found in Appendix \ref{AppE}.

\text{ }\indent  For the purposes of this study, we assume that the variance parameter $\sigma_i$ affecting the tactical uncertainty is the same for all flights and set $\sigma_i=\sigma$ for all $i\in\mathcal F$, but we consider different possible cases for $\sigma$.  For each particular value of $\sigma$, we then set the values of $\alpha_i$ and $\beta_i$ as specified in (\ref{alphabeta}) in order to ensure that the values of $(\Delta_i^{\text{pre}}+\Delta_i^{\text{tac}})$ are consistent with the historical data.  Thus, each flight has its own unique distribution for its unconstrained landing time $A_i$, and we consider different possible cases for the proportion of $\text{Var}(A_i)$ accounted for at the tactical stage as opposed to the pre-tactical stage.

\text{ }\indent We regard $\sigma$ as a primary variable of interest in this study and divide our numerical experiments into 5 strata, with the 1st, 2nd, 3rd, 4th and 5th strata consisting of experiments in which $\sigma=0.1,\;0.3,\;0.5,\;0.7,\;0.9$ respectively.  Within each stratum we randomly generate values of the other physical parameters, but the random number seeds are replicated across strata in such a way that the $n^{\text{th}}$ scenario in stratum 1 has exactly the same values of $k$, $t_0$, $t_1$, etc. as the $n^{\text{th}}$ scenario in strata 2, 3, 4 and 5 and also includes the same pre-tactical delays, flight trajectories and other scenario-specific information; in other words, these five scenarios differ only in the value of $\sigma$.  The total number of scenarios considered within each stratum is 1000, so we have 5000 scenarios in total.  We note that the smaller values of $\sigma$ result in fairly predictable ETA trajectories; for example, $\sigma=0.1$ represents the case where the unconstrained landing time of a flight expected to land in one hour has a standard deviation of about 3/4 of a minute, which could be reasonable if ETAs are not expected to vary much over time.  However, we also consider it useful to include larger $\sigma$ values in our study in order to investigate how the performances of our algorithms are affected as the amount of variability increases.

\text{ }\indent In each test scenario we sample the values of the remaining physical parameters as detailed below.  (Note: in the following list, the word `sampled' implies `sampled uniformly at random'.)\\
\begin{itemize}
\item The value of $k$ is sampled from the set $\{16,\;25,\;44,\;100,\;400\}$.  This implies that for each pair of consecutively-landing flights $(i,j)$, the separation time $M_{ij}$ has a coefficient of variation sampled from the set $\{0.05,\;0.1,\;0.15,\;0.2,\;0.25\}$.
\item With probability 0.75, the initial forecast $[t_0,\;t_1]$ for the period of bad weather is sampled from the set $\{[285,\;315],\;[270,\;330],\;[240,360]\}$.  With the remaining probability 0.25, there is no period of bad weather.
\item For simplicity, the weather variance parameter $\nu$ is set equal to $\sigma$ in all experiments, so that the uncertainty in the weather forecast is always similar to the uncertainty affecting unconstrained aircraft arrival times.
\item We set $h_i=d_i-15$ for all $i\in\mathcal F$, so that each flight enters its `tactical uncertainty' phase 15 minutes prior to its scheduled departure time.
\item For simplicity, we set $\tau=30$ in all experiments.
\item The relative cost parameter $g_i$ depends on flight $i$'s weight class in the following way: if $i$ is in the `H' class, $g_i$ is sampled from the continuous interval $[0.8,\;1]$.  If $i$ is in the `UM' or `LM' class, $g_i$ is sampled from $[0.6,\;0.8]$, and if $i$ is in the `S' class, $g_i$ is sampled from $[0.4,\;0.6]$.  This is based on the principle that heavier aircraft are likely to be carrying more passengers, and therefore should incur greater penalties for lateness.
\item The tolerance parameter $\gamma^{[S]}$ is sampled from the set $\{0,\;15\}$, with $15$ being an `industry standard' value, as mentioned in Section \ref{sec24}.  On the other hand, we simply set $\gamma^{[W]}=0$, so that any amount of air holding delay is penalized.
\item The objective function weight $\theta^{[S]}$ is sampled from the set $\{0.1,\;0.3,\;0.5,\;0.7,\;0.9\}$.  We then set $\theta^{[W]}=1-\theta^{[S]}$.\\
\end{itemize}

\text{ }\indent In the remainder of this section we compare the performances of the algorithms presented in Section \ref{sec:methodology} and also look for interesting comparisons between the characteristics of the solutions generated by these algorithms.

\subsection{Performance comparisons}\label{sec54}

\text{ }\indent Detailed performance comparisons between the various policies considered in our computational study are presented in the next three subsections.  We begin by presenting some overall summative findings, obtained by collating the results from all 5000 scenarios in our study:\\
\begin{itemize}
\item The SimHeur algorithm performed best overall, with a mean percentage improvement (with respect to the objective function (\ref{objective}), measured at the $95\%$ confidence level) of $(10.83\pm 0.41)\%$ over DetHeur and corresponding improvements of $(43.17\pm 0.35)\%$ over FCFS and $(64.99\pm 0.73)\%$ over DStat.  For clarity, the percentage improvement in a particular instance is calculated as $100\times (\Phi_{\text{H}}-\Phi_{\text{SH}})/\Phi_{\text{H}}$, where $\Phi_{\text{SH}}$ is the objective function value under SimHeur and $\Phi_{\text{H}}$ is the corresponding value under another policy $H\in\{\text{DetHeur, FCFS, DStat}\}$.
\item In 3726 of the 5000 scenarios (about 74.5\%), SimHeur's performance was better than those of DetHeur, FCFS and DStat.  In 1226 of the 5000 scenarios (about 24.5\%), DetHeur achieved the best performance.  FCFS performed best in 10 scenarios (0.2\%), and DStat performed best in 38 scenarios (about 0.8\%)\\
\end{itemize}

\text{ }\indent Our objective function (\ref{objective}) is of a quadratic nature, and also depends on the scales used for the weights $\theta^{[S]}$ and $\theta^{[W]}$ and the relative cost parameters $g_i$ for $i\in\mathcal F$, so its values are not particularly easy to interpret in a physical context.  Therefore, in order to indicate how SimHeur compares to DetHeur, FCFS and DStat with respect to more interpretable performance measures related to schedule punctuality and flight duration, we provide some additional summary statistics below.  (All of the reported confidence intervals are at the $95\%$ level.)\\
\begin{itemize}
\item The mean schedule delay, defined as the difference between the actual landing time $L_i$ and the scheduled landing time $a_i$ averaged over all flights $i\in\mathcal F$ and all 5000 scenarios, was $(6.19\pm 0.03)$ minutes for the SimHeur algorithm.  For the DetHeur, FCFS and DStat algorithms, the corresponding results were $(7.64\pm 0.03)$, $(10.10\pm 0.04)$ and $(25.81\pm 0.04)$ minutes, respectively.
\item The mean airborne holding delay, defined as the difference between the actual landing time $L_i$ and the landing time $A_i-\rho_i$ that would occur if flight $i$ was released from the pool immediately and avoided any queueing delay (averaged over all flights $i\in\mathcal F$ and all 5000 scenarios), was $(15.28\pm 0.02)$ minutes for the SimHeur algorithm.  For the DetHeur, FCFS and DStat algorithms, the corresponding results were $(16.73\pm 0.02)$, $(19.19\pm 0.02)$ and $(34.90\pm 0.04)$ minutes, respectively.  As explained in Section \ref{sec24}, this holding delay can be incurred either in the pool or in the queue (or both).  In the case of the FCFS algorithm, it is only incurred in the queue.\\
\end{itemize}

Unsurprisingly, the relative performances of these various policies are influenced strongly by the values of the physical parameters.  In the remainder of this subsection we provide tables to show how the comparisons are affected by varying these parameters.\\

\subsubsection{Varying $\sigma$}~
\text{ }\indent Table \ref{tab_sigma} shows how the comparisons are affected by adjusting the value of the tactical variance parameter $\sigma$.  Recall that we also set $\nu=\sigma$ in all experiments.  In Table \ref{tab_sigma} (and also Tables \ref{tab_k}-\ref{tab_theta} later), columns 3, 4 and 5 show 95\% confidence intervals for the mean percentage improvement (with respect to the value of the objective function (\ref{objective})) achieved by SimHeur against DetHeur, FCFS and DStat respectively, while columns 6 (resp. 7, 8, 9) show the percentage of all experiments in which SimHeur (resp. DetHeur, FCFS, DStat) achieved the smallest objective function value.

{\renewcommand{\arraystretch}{1.25}%
\begin{table}[htbp]
  \centering \small
  \caption{ Comparisons between SimHeur (SH), DetHeur (DH), first-come-first-served (FCFS) and DStat respectively, for various values of the tactical variance parameter $\sigma$.}
    \begin{tabular}{|c|c|c|c|c|c|c|c|c|}
\cline{3-9}    \multicolumn{1}{c}{} &       & \multicolumn{3}{c|}{\textbf{Pct. Improvement}} & \multicolumn{4}{c|}{\textbf{Pct. of Experiments}} \bigstrut\\
    \hline
    $\sigma$ value & Count & \footnotesize SH vs. DH & \footnotesize SH vs. FCFS & \footnotesize SH vs. DStat & \footnotesize SH best & \footnotesize DH best & \footnotesize FCFS best & \footnotesize DStat best\bigstrut\\
    \hline
    $0.1$ & 1000   & $2.14\pm 0.56$ & $51.20\pm 0.62$ & $19.70\pm 0.78$ & 56.30 & 39.90 & 0.00 & 3.80\\
    \rowcolor{Grey}
    $0.3$ & 1000   & $7.56\pm 0.72$ & $48.33\pm 0.63$ & $59.54\pm 0.76$ & 75.30 & 24.70 & 0.00 & 0.00\\
    $0.5$ & 1000   & $12.43\pm 0.88$ & $44.49\pm 0.64$ & $76.22\pm 0.53$ & 79.50 & 20.50 & 0.00 & 0.00\\
    \rowcolor{Grey}
    $0.7$ & 1000   & $15.22\pm 0.96$ & $38.95\pm 0.70$ & $83.19\pm 0.41$ & 80.10 & 19.80 & 0.10 & 0.00\\
    $0.9$ & 1000   & $16.82\pm 1.01$ & $32.85\pm 0.75$ & $86.27\pm 0.33$ & 81.40 & 17.70 & 0.90 & 0.00\bigstrut[b]\\
    \hline
    \end{tabular}%
  \label{tab_sigma}%
\end{table}%

\text{ }\indent One would expect the relative improvement of SimHeur versus DetHeur to increase as the amount of stochasticity in our model increases.  Indeed, Table \ref{tab_sigma} shows a clear trend for SimHeur to improve its advantage over DetHeur as $\sigma$ is increased.  On the other hand, SimHeur's advantage over FCFS diminishes when $\sigma$ is increased.  This seems to indicate that the FCFS policy becomes stronger when there is more uncertainty in the unconstrained arrival times $A_i$, which can be explained by the fact that flights are more likely to arrive significantly later than their scheduled landing times $a_i$ in such circumstances, and must then be released as soon as possible in order to avoid large penalties.  The DStat algorithm performs relatively well in the $\sigma=0.1$ case and it can even outperform the other policies in some instances, which is a sign that the heuristic procedure described in Appendix \ref{AppD} works well for solving the relevant static, deterministic optimization problem.  However, its performance deteriorates sharply as $\sigma$ increases and it becomes worse than FCFS in such cases.  Indeed, larger $\sigma$ values can cause DStat to incur very large holding costs due to the increased likelihood of flights arriving in the pool either very late (and thus causing other flights to wait before being released) or very early (and thus having to wait themselves).

\subsubsection{Varying $k$}~

Table \ref{tab_k} shows how the comparisons are affected by adjusting the value of the variance parameter $k$ for the Erlang-distributed separation times.

{\renewcommand{\arraystretch}{1.25}%
\begin{table}[htbp]
  \centering \small
  \caption{ Comparisons between SimHeur (SH), DetHeur (DH), first-come-first-served (FCFS) and DStat for various values of the service time variance parameter $k$.}
    \begin{tabular}{|c|c|c|c|c|c|c|c|c|}
\cline{3-9}    \multicolumn{1}{c}{} &       & \multicolumn{3}{c|}{\textbf{Pct. Improvement}} & \multicolumn{4}{c|}{\textbf{Pct. of Experiments}} \bigstrut\\
    \hline
    $k$ value & \footnotesize Count & \footnotesize SH vs. DH & \footnotesize SH vs. FCFS & \footnotesize SH vs. DStat & \footnotesize SH best & \footnotesize DH best & \footnotesize FCFS best & \footnotesize DStat best \bigstrut\\
    \hline
    $16$ & 1000   & $11.80\pm 0.92$ & $42.17\pm 0.80$ & $64.40\pm 1.62$ & 75.50 & 23.70 & 0.40 & 0.40 \\
    \rowcolor{Grey}
    $25$ & 972   & $11.52\pm 0.96$ & $42.20\pm 0.82$ & $64.75\pm 1.69$ & 73.35 & 25.82 & 0.21 & 0.62 \\
    $44$ & 1039   & $11.66\pm 0.98$ & $43.21\pm 0.76$ & $64.76\pm 1.63$ & 71.70 & 26.95 & 0.29 & 1.06 \\
    \rowcolor{Grey}
    $100$ & 969   & $9.60\pm 0.84$ & $43.83\pm 0.78$ & $65.14\pm 1.66$ & 75.75 & 22.91 & 0.00 & 1.34 \\
    $400$ & 1020   & $9.56\pm 0.80$ & $44.38\pm 0.75$ & $65.87\pm 1.53$ & 76.37 & 23.14 & 0.10 & 0.39 \bigstrut[b]\\
    \hline
    \end{tabular}%
  \label{tab_k}%
\end{table}%

\text{ }\indent We recall that smaller $k$ values are associated with more stochasticity in service times and, although there is a trend for SimHeur's advantage over DetHeur to diminish as $k$ increases, this trend is much less significant than the one observed in Table \ref{tab_sigma} (for the $\sigma$ values).  The FCFS policy seems to perform slightly better (relative to SimHeur) when $k$ is smaller; indeed, it intuitively makes sense that we would want to take a conservative approach and release flights early if their service times are subject to a lot of uncertainty. 

\subsubsection{Varying $[t_0,\;t_1]$}~

Table \ref{tab_t0_t1} shows how the comparisons are affected by adjusting the initial forecast, $[t_0,\;t_1]$, for the period of bad weather.  Recall that the actual duration of the bad weather period is subject to random variation, as described in Section \ref{sec24}.

{\renewcommand{\arraystretch}{1.25}%
\begin{table}[htbp]
  \centering \small
  \caption{ Comparisons between SimHeur (SH), DetHeur (DH), first-come-first-served (FCFS) and DStat for various different intervals $[t_0,t_1]$ (note that the case $\emptyset$ represents `no bad weather').}
    \begin{tabular}{|c|c|c|c|c|c|c|c|c|}
\cline{3-9}    \multicolumn{1}{c}{} &       & \multicolumn{3}{c|}{\textbf{Pct. Improvement}} & \multicolumn{4}{c|}{\textbf{Pct. of Experiments}} \bigstrut\\
    \hline
    $[t_0,t_1]$ interval & \footnotesize Count & \footnotesize SH vs. DH & \footnotesize SH vs. FCFS & \footnotesize SH vs. DStat & \footnotesize SH best & \footnotesize DH best & \footnotesize FCFS best & \footnotesize DStat best \bigstrut\\
    \hline
    $\emptyset$ & 1286   & $10.75\pm 0.79$ & $39.17\pm 0.71$ & $66.93\pm 1.41$ & 74.65 & 24.49 & 0.54 & 0.31 \\
    \rowcolor{Grey}
    $[285,315]$ & 1228   & $10.74\pm 0.83$ & $42.14\pm 0.73$ & $66.65\pm 1.44$ & 74.67 & 25.00 & 0.16 & 0.16 \\
    $[270,330]$ & 1232   & $10.69\pm 0.81$ & $44.44\pm 0.67$ & $64.51\pm 1.48$ & 74.03 & 24.84 & 0.08 & 1.06 \\
    \rowcolor{Grey}
    $[240,360]$ & 1254   & $11.14\pm 0.81$ & $47.02\pm 0.61$ & $61.82\pm 1.47$ & 74.72 & 23.76 & 0.00 & 1.52 \bigstrut[b]\\
    \hline
    \end{tabular}%
  \label{tab_t0_t1}%
\end{table}%

\text{ }\indent In this case there is no obvious trend for SimHeur's advantage over DetHeur to become larger or smaller as the duration of the bad weather period (according to the initial forecast) is increased.  Instead, it robustly maintains an advantage of about $11\%$ (on average) over all scenarios.  The FCFS policy seems to become worse as the duration of bad weather increases, which can be explained by the fact that in bad weather, it may be advantageous to delay a flight's release from the pool due to the dependence of the required separation time on weather conditions at the time of release.  The same principle may be used to explain the DStat policy's apparent improvement as the duration of bad weather increases, as this policy has a tendency to impose long pool-holding times which might help some flights to avoid landing during bad weather. 

\subsubsection{Varying $\gamma^{[S]}$}~

Table \ref{tab_gamma} shows how the comparisons are affected by adjusting the `tolerance' parameter $\gamma^{[S]}$ in the objective function (\ref{objective}).  For this parameter we only tested the cases $\gamma^{[S]}=0$ and $\gamma^{[S]}=15$ (measured in minutes).

\begin{table}[htbp]
  \centering
  \caption{ Comparisons between SimHeur (SH), DetHeur (DH), first-come-first-served (FCFS) and DStat for the two cases $\gamma^{[S]}=0$ and $\gamma^{[S]}=15$.}
    \begin{tabular}{|c|c|c|c|c|c|c|c|c|}
\cline{3-9}    \multicolumn{1}{c}{} &       & \multicolumn{3}{c|}{\textbf{Pct. Improvement}} & \multicolumn{4}{c|}{\textbf{Pct. of Experiments}} \bigstrut\\
    \hline
    $\gamma^{[S]}$ value & \footnotesize Count & \footnotesize SH vs. DH & \footnotesize SH vs. FCFS & \footnotesize SH vs. DStat & \footnotesize SH best & \footnotesize DH best & \footnotesize FCFS best & \footnotesize DStat best \bigstrut\\
    \hline
    $0$ & 2525   & $12.30\pm 0.56$ & $40.47\pm 0.48$ & $63.62 \pm 1.02$ & 77.98 & 20.99 & 0.28 & 0.75 \\
    \rowcolor{Grey}
    $15$ & 2475   & $9.34\pm 0.58$ & $45.91\pm 0.49$ & $66.38\pm 1.03$ & 70.99 & 28.12 & 0.12 & 0.77 \bigstrut[b]\\
    \hline
    \end{tabular}%
  \label{tab_gamma}%
\end{table}%

\text{ }\indent From Table \ref{tab_gamma} we infer that SimHeur's advantage over DetHeur is greater in the case $\gamma^{[S]}=0$.  In this case it becomes more critical for flights to land near their scheduled times, which tends to imply that a FCFS policy will perform better (as suggested by the results in the 4th column).  The DetHeur policy tends to be overly optimistic about the actual time required for a sequence of flights to complete service, so its weaker performance in the $\gamma^{[S]}=0$ case may be explained by its tendency to keep flights waiting for too long in the pool.  

\subsubsection{Varying $\theta^{[S]}$ and $\theta^{[W]}$}~

Table \ref{tab_theta} shows how the comparisons are affected by adjusting the objective function weights $\theta^{[S]}$ and $\theta^{[W]}$, subject to the constraint $\theta^{[S]}+\theta^{[W]}=1$.

{\renewcommand{\arraystretch}{1.25}%
\begin{table}[htbp]
  \centering \small
  \caption{ Comparisons between SimHeur (SH), DetHeur (DH), first-come-first-served (FCFS) and DStat for various values of the objective function weights $\theta^{[S]}$ and $\theta^{[W]}=1-\theta^{[S]}$.}
    \begin{tabular}{|c|c|c|c|c|c|c|c|c|}
\cline{3-8}    \multicolumn{1}{c}{} &       & \multicolumn{3}{c|}{\textbf{Pct. Improvement}} & \multicolumn{4}{c|}{\textbf{Pct. of Experiments}} \bigstrut\\
    \hline
    $\theta^{[S]}$ value & \footnotesize Count & \footnotesize SH vs. DH & \footnotesize SH vs. FCFS & \footnotesize SH vs. DStat & \footnotesize SH best & \footnotesize DH best & \footnotesize FCFS best & \footnotesize DStat best \bigstrut\\
    \hline
    $0.1$ & 1015   & $-2.59\pm 0.62$ & $45.09\pm 0.72$ & $67.89\pm 1.62$ & 36.45 & 62.86 & 0.00 & 0.69 \\
    \rowcolor{Grey}
    $0.3$ & 995   & $7.96\pm 0.80$ & $43.62\pm 0.71$ & $67.46\pm 1.56$ & 70.95 & 28.54 & 0.00 & 0.50 \\
    $0.5$ & 995   & $15.71\pm 0.85$ & $42.55\pm 0.74$ & $66.21\pm 1.56$ & 87.64 & 11.86 & 0.30 & 0.20 \\
    \rowcolor{Grey}
    $0.7$ & 983   & $17.63\pm 0.84$ & $42.81\pm 0.81$ & $63.11\pm 1.64$ & 91.15 & 7.73 & 0.51 & 0.61 \\
    $0.9$ & 1012   & $15.73\pm 0.77$ & $41.74\pm 0.91$ & $60.26\pm 1.69$ & 87.15 & 10.87 & 0.20 & 1.78 \bigstrut[b]\\
    \hline
    \end{tabular}%
  \label{tab_theta}%
\end{table}%

\text{ }\indent In Table \ref{tab_theta} we observe (for the first time in our numerical study) a subset of experiments in which SimHeur performs worse than DetHeur.  Specifically, when $\theta^{[S]}=0.1$ and $\theta^{[W]}=0.9$, SimHeur appears to be relatively weak.  On the other hand, it becomes very strong (relative to DetHeur) when $\theta^{[S]}$ is increased.  In order to find a possible explanation for this, we recall that SimHeur's advantage over DetHeur is usually derived from the fact that it evaluates the costs of potential runway sequences more accurately than DetHeur under stochastic conditions, whereas DetHeur has the advantage of being able to explore the solution space faster than SimHeur (as discussed in Section \ref{sec32}).  When $\theta^{[S]}$ is small relative to $\theta^{[W]}$, it becomes less important for flights to land near their scheduled times and the problem is mainly about arranging the sequence of landings so that the total time required for all services to complete (referred to as a `makespan' in the job shop literature) is minimized.  This depends a lot on controlling the sequence of aircraft weight classes in order to reduce average separation times, which implies that a FCFS sequence tends to become weaker (as can be seen from the 4th column in Table \ref{tab_theta}), as it pays no attention to these weight classes.  Our conjecture is that in the $\theta^{[S]}=0.1$ case, the sequencing problem becomes dominated by the aircraft weight classes rather than the on-time requirements for individual flights, and the optimal sequence tends to require larger deviations from a FCFS sequence.  Under such circumstances, DetHeur is able to outperform SimHeur because it explores the solution space faster and performs many more mutation steps, enabling it to discover better sequences.

\subsection{Sequencing patterns with respect to aircraft weight classes}\label{sec55}

\text{ }\indent It is well-known in aircraft sequencing problems that in order to minimize average separation times (or, equivalently, maximize average throughput rates) one should aim to have `strings' of the same aircraft type appearing successively in the landing sequence; for example, the sequence `H-H-H-LM-LM-LM' will have a smaller expected average separation time than the sequence `H-LM-H-LM-H-LM' if one assumes that these aircraft are within sufficient proximity of each other to aim for the smallest allowable time separations (this can be verified using the required separation times in Table \ref{sep_times}).  Of course, the FCFS algorithm will ignore this principle, but one might expect the other heuristics to implement sequences in which leader-follower pairs are of the same weight class more frequently than one would expect under a FCFS sequence.

\text{ }\indent Table \ref{tab_follow_pct} shows a comparison between the different algorithms with respect to the average percentage of flights in $\mathcal F$ that immediately follow another flight of the same weight class in the landing sequence (where the average is taken over all 5000 scenarios).

{\renewcommand{\arraystretch}{1.25}%
\begin{table}[htbp]
  \centering \small
  \caption{ Comparisons between SimHeur, DetHeur, first-come-first-served (FCFS) and DStat with respect to the percentage of flights landing immediately after another flight of the same weight class, for various values of $\sigma$.}
    \begin{tabular}{|c|c|c|c|c|c|}
    \hline
    $\sigma$ value & Count & SimHeur & DetHeur & FCFS & DStat \bigstrut\\
    \hline
    $0.1$ & 1000   & $75.29$ & $75.99$ & $56.52$ & 79.12 \\
    \rowcolor{Grey}
    $0.3$ & 1000   & $75.53$ & $76.56$ & $56.62$ & 79.30 \\
    $0.5$ & 1000   & $75.91$ & $77.44$ & $56.57$ & 79.79 \\
    \rowcolor{Grey}
    $0.7$ & 1000   & $75.94$ & $78.04$ & $56.76$ & 80.82 \\
    $0.9$ & 1000   & $75.65$ & $78.48$ & $56.62$ & 81.55  \bigstrut[b]\\
    \hline
    \end{tabular}%
  \label{tab_follow_pct}%
\end{table}%

\text{ }\indent It is clear from Table \ref{tab_follow_pct} that the SimHeur and DetHeur algorithms are able to derive an advantage over the FCFS algorithm by `grouping' flights of the same weight class together, in the manner described above.  The DetHeur algorithm tends to do this to a slightly greater extent than SimHeur, which may be explained by its tendency to underestimate the costs of potential runway sequences (and hence overestimate the amount of available time to arrange flights into suitable `strings' before they need to be released from the pool).  A similar explanation can be given for the fact that the percentages for DStat are even higher than those for SimHeur and DetHeur, as the DStat algorithm is based on making predictions over the entire interval $\mathcal T$, without any recourse to actual observed events; thus, it will tend to be even more optimistic than DetHeur with regard to tactical delays and the amount of available time for arranging the weight class sequence.

\text{ }\indent The value of $\sigma$ has relatively little effect on these findings, although it is interesting to note that the percentages for DStat increase as $\sigma$ increases, despite the fact that DStat does not change its runway sequence in response to any random events during $\mathcal T$.  Recall that we assume DStat has knowledge of the pre-tactical delays $\Delta_i^{\text{pre}}$ (realized at the beginning of $\mathcal T$) and, as $\sigma$ becomes larger, the pre-tactical delays become more predictable because we suppress the amount of pre-tactical uncertainty (as opposed to tactical uncertainty).  The arrivals schedule (shown in Figure \ref{sec4_fig}) shows that many `heavy' flights are due to arrive consecutively in the early morning hours.  For smaller values of $\sigma$, DStat's predictions are based on a greater amount of `mixing' between the expected positions in the arrival sequence of different weight classes due to increased pre-tactical uncertainty.  Hence, DStat accounts for this in its sequencing decisions and perceives less opportunity to group flights of the same weight class together. 

\subsection{Average times spent in the pool under SimHeur and DetHeur}\label{sec56}

\text{ }\indent We can also compare the various solution algorithms with respect to the average time that a flight spends in the pool before being released and added to the runway queue.  The entries in Table \ref{tab_av_hold} are obtained by averaging over all flights $i\in\mathcal F$ and all 5000 scenarios. 
We do not include the FCFS algorithm in the table as it obviously does not keep any flights waiting in the pool.

{\renewcommand{\arraystretch}{1.25}%
\begin{table}[htbp]
  \centering \small
  \caption{ Comparisons between SimHeur, DetHeur and DStat with respect to the average holding time in the pool (in minutes), for various values of $\sigma$.}
    \begin{tabular}{|c|c|c|c|c|c|}
    \hline
    $\sigma$ value & Count & SimHeur & DetHeur & DStat \bigstrut\\
    \hline
    $0.1$ & 1000   & $6.77\pm 0.03$ & $3.39\pm 0.02$ & $3.69\pm 0.02$  \\
    \rowcolor{Grey}
    $0.3$ & 1000   & $6.02\pm 0.03$ & $3.99\pm 0.02$ & $7.16\pm 0.03$ \\
    $0.5$ & 1000   & $5.34\pm 0.03$ & $4.51\pm 0.02$ & $12.26\pm 0.04$ \\
    \rowcolor{Grey}
    $0.7$ & 1000   & $4.88\pm 0.03$ & $5.04\pm 0.02$ & $18.49\pm 0.05$ \\
    $0.9$ & 1000   & $4.70\pm 0.03$ & $5.79\pm 0.03$ & $25.47\pm 0.07$ \bigstrut[b]\\
    \hline
    \end{tabular}%
  \label{tab_av_hold}%
\end{table}%

\text{ }\indent There are several interesting features to discuss in Table \ref{tab_av_hold}.  It appears that SimHeur tends to release flights from the pool earlier as $\sigma$ increases, but the opposite is true for DetHeur.  In the case of SimHeur, we recall that flights are not released from the pool unless the condition $\Upsilon_{s^*}^{(n)}>\lambda$ is satisfied (see Section \ref{sec31}).  This essentially means that we delay the release until the estimated probability of the released plane going `straight into service' (i.e. avoiding a queueing delay) exceeds the threshold $\lambda$, which we set to zero in these experiments.  Larger values of $\sigma$ imply more variability in the additional travel time $A_i-R_i$ needed to progress to the runway, and therefore the algorithm will tend to sample fewer random trajectories before discovering a particular trajectory under which the flight goes straight into service, resulting in the condition $\Upsilon_{s^*}^{(n)}>\lambda$ being met.  On the other hand, DetHeur has no comparable mechanism to delay flight releases and instead releases a flight immediately if it is already in the pool and judged to be the best one to release next.  As $\sigma$ increases, flights become more likely to enter the pool far in advance of their scheduled landing times, which could explain the tendency for DetHeur to delay them in the pool for longer.  

\text{ }\indent It is also clear that DStat becomes susceptible to very large pool-holding delays as $\sigma$ increases; this is because flights that arrive in the pool very early or very late will tend to cause long pool-holding delays either to themselves or to other flights.  Realistically, the DStat policy would probably not be implementable at all under highly stochastic conditions.

\subsection{Cost predictions for runway sequences under SimHeur and DetHeur}\label{sec57}

\text{ }\indent The SimHeur and DetHeur algorithms make decisions based on continuously evaluating the estimated costs associated with different possible runway sequences, as described in Section \ref{sec:methodology}.  Since the DetHeur algorithm is based on the assumption that all unrealized random variables will conform to their expected values, its predictions for runway sequence costs are likely to be overly optimistic.  On the other hand, since SimHeur's cost estimates are based on the simulation of many possible random trajectories, its predictions are likely to be much more accurate.  Using additional detailed output from our experiments (omitted from this paper, but obtainable from the authors upon request) we have been able to verify that SimHeur's cost predictions tend to be both higher and more accurate than DetHeur's - although, on rare occasions, SimHeur can produce absurdly high estimates due to one or more of its randomly-sampled trajectories being associated with extreme sequences of delays; i.e. it is prone to outliers.

\text{ }\indent The fact that DetHeur's decision-making process is based on relatively inaccurate predictions of the costs obtained under possible alternative sequences raises some interesting questions.  Of course, the tendency to underestimate sequence costs does not, in itself, imply that the algorithm performs badly.  Indeed, the results in Section \ref{sec54} have already shown that DetHeur is capable of performing well in many test scenarios, particularly when $\sigma$ is small.  Clearly, the most important task for SimHeur and DetHeur in our problem is to rank sequences in the correct order, rather than to predict the costs accurately.  If both algorithms tend to rank sequences in the same order over the course of time, then they will tend to achieve similar performances, even if one is predicting costs much more accurately than the other.  However, we might consider ways to modify the problem in such a way that predictions of sequence costs are linked to other decision options.  For example, suppose there is an option to somehow make `extra capacity' available for a limited period of time $-$ by using an extra runway for arrivals, for example $-$ or to reduce demand rates by diverting or canceling some flights, subject to penalty costs.  In such circumstances, the SimHeur algorithm may be much better than DetHeur in assessing when such congestion-mitigating actions are necessary and we might see a wider performance gap between the algorithms, even when $\sigma$ is small.  We leave this as a possible direction for further research.\\

\section{Conclusions}\label{sec:conclusions}

\text{ }\indent This paper has introduced an original mathematical model for stochastic runway scheduling, in which aircraft ETAs and weather conditions evolve dynamically according to continuous-time stochastic processes, while runway `service times' depend on sequence-dependent Erlang distributions.  The aim is to consider a high-dimensional and information-rich environment in which air traffic controllers are able to update their plans frequently in response to the latest unfolding events.  It is natural to consider a simheuristic approach to such a problem, since other conventional optimization approaches (e.g. two-stage stochastic programming) are not well-equipped to deal with the continuous nature of the information updates and decision epochs in our model.

\text{ }\indent Our numerical experiments, based on a schedule of more than 300 arrivals at Heathrow Airport and configured using a large set of historical on-time performance data, have shown that the proposed simheuristic algorithm (SimHeur) is capable of outperforming an alternative method based on deterministic forecasts (DetHeur) under a wide range of parameter values, and also improves substantially upon a simple `first-come-first-served' policy.  Notably, even when the amount of stochasticity in our model is relatively low (e.g. with the choices $\sigma=\nu=0.1$ and $k=400$), we find that the improvements given by SimHeur versus DetHeur are significant, and these improvements would likely be greater without the severe time compression used in our experiments (i.e. one minute of CPU time to represent one hour of real time).  It should also be noted that the advantage of SimHeur over DetHeur tends to become greater when the on-time requirements of individual flights are given more weight in the objective function.

\text{ }\indent Certainly, our computational study could be expanded in order to consider other schedules, alternative objective functions and different model parameters (including both `physical' parameters related to real-world operating conditions and the tuning parameters used by our algorithms), and this should be a direction of further research.  We also plan to consider models with both arrivals and departures in future work, and to consider decision-making problems in which accurate predictions of future costs can help to ensure that extra airport resources or demand management strategies are deployed at the most critical times, as described in Section \ref{sec57}.\\

\textbf{Acknowledgments.  }This work has been supported by the Engineering and Physical Sciences Research Council (EPSRC) through Programme Grant EP/M020258/1 ``Mathematical models and algorithms for allocating scarce airport resources (OR-MASTER)".  We would also like to thank Flightradar24 (https://www.flightradar24.com/) for granting permission for their on-time performance data to be used in our study.\\

\bibliographystyle{apalike} 
\bibliography{Simheuristics_Arxiv_paper} 

\newpage
\Large \textbf{Appendices}
\normalsize

\appendix
\normalsize

\normalsize
 \section{The heuristic move operator $H$}\label{AppA}

\text{ }\indent Given a sequence $s$ of length $l$, the heuristic move operator $H$ works as follows:\\
\begin{enumerate}
\item Let $L:=l(l+1)/2$.  We randomly select an integer $j\in\{1,2,...,l\}$ in the following way: the probability of selecting 1 is $l/L$, the probability of selecting 2 is $(l-1)/L$, and in general, the probability of selecting $p$ is $(l+1-p)/L$ for $1\leq p\leq l$.
\item The flight in position $j$ of $s$ is to be shifted by a certain number of positions either forwards or backwards.  First, decide on the direction of movement using a simple `coin flip', so that the `forwards' and `backwards' directions are selected with probability $0.5$ each.
\item Let $P$ be sampled uniformly at random from the set $\{1,2,3\}$.  The flight in position $j$ of $s$ is removed from $s$ and then re-inserted at position $\max\{j-P,\;1\}$ (if the `forwards' direction was selected in step 2) or position $\min\{j+P,\;l\}$ (if the `backwards' direction was selected).\\
\end{enumerate}

\text{ }\indent We note that if the integer $j$ selected in step 1 of the above procedure is either $1$ or $l$, then there is a possibility that the `new' sequence generated is the same as the old one.  However, the nature of how $H$ is used in the SimHeur and DetHeur algorithms implies that the procedure will be repeated as many times as necessary until the required number of distinct sequences have been found.\\

 \section{The method of pre-generating random events}\label{AppB}

 \text{ }\indent The pre-generated sample path $\Omega^{\text{true}}$ includes complete trajectories $\{X_i(t)\}_{t\in\mathcal T}$ for $i\in\mathcal F$, weather forecast trajectories $\{T_0(t)\}_{t\in\mathcal T}$ and $\{T_1(t)\}_{t\in\mathcal T}$, and a set of possible separation times $\{M_{ij}\}$ for each $j\in\mathcal F$ (taking into account all possibilities for the weight class of the predecessor flight $i$).  The approach for generating the trajectories $\{X_i(t)\}_{t\in\mathcal T}$ relies upon a time discretization.  The steps are described below.\\

\begin{enumerate}
\item For each $i\in\mathcal F$ we set the initial ETA, $X_i(0)$, as follows:
$$X_i(0)=\begin{cases}a_i+\Delta_i^{\text{pre}}+\text{N}(0,\;\sigma_i^2 |h_i|),&\text{if }h_i<0,\\
a_i+\Delta_i^{\text{pre}},&\text{if }h_i\geq 0,\nonumber\end{cases}$$
where $\text{N}(a,b)$ denotes a Normal random variable with mean $a$ and variance $b$.
\item We step forward in hundredths of a minute and, for each $t=0.01,\;0.02,\;...,T$, generate the value $X_i(t)$ using a standard method for simulating Brownian motion:
$$X_i(t):=X_i(t-0.01)+\text{N}(0,\;0.01\times\sigma_i^2).$$
We also define $X_i(u):=X_i(t-0.01)$ for $u\in[t-0.01,t)$.\\
\end{enumerate}

\text{ }\indent The above method ensures that $X_i(t)$ is defined for all $t\in\mathcal T$ (it is a piecewise constant function of time).  Although a time discretization is used, the SimHeur and DetHeur algorithms may carry out system state updates at any time $t$ in the continuous interval $\mathcal T$ (with $t$ depending on the exact amount of CPU time spent on computations).  In practice, this means that the values $X_i(0)$, $X_i(0.01)$, $X_i(0.02)$, etc. are stored inside an array and at time $t$, the algorithm looks up the value $X_i(t')$ where $t':=\max\{u\in\{0,\;0.01,\;0.02,...,T\}:u\leq t\}$; for example, if $t=2.764$ then we look up the value $X_i(2.76)$.

\text{ }\indent The steps for generating $\{T_0(t)\}_{t\in\mathcal T}$ and $\{T_1(t)\}_{t\in\mathcal T}$ are very similar to the above, except we set initial values $T_0(0)=t_0$ and $T_1(0)=t_1$ and use $\nu$ as the variance parameter instead of $\sigma_i$.  

\text{ }\indent In order to simulate separation times $M_{ij}$, we pre-generate a set of values $\{p_j^{\text{sep}}\}_{j\in\mathcal F}$, with each $p_j^{\text{sep}}$ being randomly sampled from a $\text{Uniform}[0,1]$ distribution.  Then, once the weight class of the preceding flight $i$ becomes known during the simulation, we calculate $M_{ij}$ by sampling the $(p_j^{\text{sep}})^{\text{th}}$ quantile of the Erlang distribution for $M_{ij}$.\\

 \section{The mutation step}\label{AppC}

\text{ }\indent  Given a sequence $s$ of length $l$, the steps used by the SimHeur and DetHeur algorithms to perform a `mutation' of this sequence at a particular time $t\in\mathcal T$ are as follows:\\

\begin{enumerate}
\item Let $M(t)$ denote the number of flights in $\mathcal F$ that have not yet been added to the landing queue at time $t$; that is, $M(t)$ consists of flights that are either still in the pool or yet to arrive in the pool.  
\item Let $P(t):=\min\{4,\;Q(t)\}$, where $Q(t):=\min\{l,\;M(t)\}$.  Here, $P(t)$ is interpreted as the length of a particular subsequence within $s$ that we want to `shuffle' in order to obtain a new sequence.
\item Define $R(t):=Q(t)-P(t)+1$ as the number of possible starting positions of the string that we are going to shuffle.
\item Let $L:=R(t)(R(t)+1)/2$.  We randomly select an integer $j\in\{1,2,...,R(t)\}$ in the following way: the probability of selecting 1 is $R(t)/L$, the probability of selecting 2 is $(R(t)-1)/L$, and in general, the probability of selecting $p$ is $(R(t)+1-p)/L$ for $1\leq p\leq R(t)$.
\item Consider the subsequence formed by taking the flights in positions $j,\;j+1,...,j+P(t)-1$ of $s$.  Remove all flights in this subsequence from $s$, then `shuffle' the subsequence, i.e. choose a random permutation of it.  Finally, re-insert the shuffled subsequence in the same position within $s$ that it occupied before.\\
\end{enumerate}

\text{ }\indent Following these steps, we obtain a new sequence, interpreted as a `mutation' of $s$.\\

 \section{Obtaining the DStat policy}\label{AppD}

 \text{ }\indent The heuristic method for obtaining the DStat policy referred to in Section \ref{sec33} is as follows:\\

\begin{enumerate}
\item We begin with a sequence $s_0$ of length $|\mathcal F|$ in which the flights are ordered according to their ETAs following the realization of pre-tactical uncertainty; that is, if flight $i\in\mathcal F$ appears before flight $j\in\mathcal F$ then this implies $a_i+\Delta_i^{\text{pre}}\leq a_j+\Delta_j^{\text{pre}}$.
\item Set $s^{\text{best}}:=s_0$ and initialize a counter $c=0$.
\item Set the pool arrival time to $Q_i=a_i+\Delta_i^{\text{pre}}-\tau$ for all $i\in\mathcal F$, set $U_0=t_0$, $U_1=t_1$ and assume that all separation times $M_{ij}$ conform to their expected values $E_{ij}(R_j)$ (defined in (\ref{Eij_eq})).  For the release times $R_j$, we assume that each flight $i$ in $s^{\text{best}}$ is released from the pool at the earliest possible moment, with the restriction that flights must be released in the same order that they appear in $s^{\text{best}}$.  For example, if the 2nd flight in $s^{\text{best}}$ happens to arrive in the pool earlier than the 1st flight, then the 2nd flight should be released immediately after the 1st flight arrives in (and is released from) the pool.  We also define $A_j=R_j+\tau$ as the unconstrained arrival time of flight $j\in\mathcal F$, with the actual landing time given by $L_j=\max\{A_j,L_i+M_{ij}\}$, where $i$ is the predecessor of $j$ in the queue.  We then calculate the value of the objective function (\ref{objective}) if the sequence $s^{\text{best}}$ is followed and let this be denoted by $C^{\text{best}}$.
\item Select a subsequence of 6 flights appearing consecutively in $s^{\text{best}}$ by sampling uniformly at random from the $|\mathcal F|-5$ possible alternatives (corresponding to starting positions $1,2,...,|\mathcal F|-5$ within $s^{\text{best}}$).
\item The subsequence selected in step 4 is removed from $s^{\text{best}}$ and the flights in this subsequence are then `shuffled'; i.e. a random permutation of the subsequence is chosen.  The shuffled subsequence is then re-inserted into $s^{\text{best}}$ in the same position that it occupied before.  This yields a new sequence that we refer to as $s^{\text{curr}}$.
\item Calculate the value of the objective function (\ref{objective}) under the new sequence $s^{\text{curr}}$ by setting the values of $Q_i$, $U_0$, $U_1$, $M_{ij}$ to their expected values and calculating $R_i$, $A_i$ as described in step 3.  Let $C^{\text{curr}}$ denote the objective function value under $s^{\text{curr}}$.
\item If $C^{\text{curr}}<C^{\text{best}}$, then set $C^{\text{best}}:=C^{\text{curr}}$ and $s^{\text{best}}:=s^{\text{curr}}$, and set $c=0$.  Otherwise, increase the counter $c$ by 1.
\item If $c=10,000$, then terminate the procedure and accept $s^{\text{best}}$ as the DStat policy.  Otherwise, return to step 4.\\
\end{enumerate}

 \text{ }\indent Note that the final value of $C^{\text{best}}$ associated with the sequence $s^{\text{best}}$ is not a measure of DStat's actual performance, as $C^{\text{best}}$ is based on the assumption of all random variables conforming to their expected values.  To evaluate DStat's actual performance, we must use the `real' information contained in $\Omega^{\text{true}}$, as explained in Section \ref{sec33}.\\
 
  \section{Values of algorithmic parameters}\label{AppE}

 \text{ }\indent In order to determine suitable values for the algorithmic parameters used by the SimHeur and DetHeur algorithms, we conducted a preliminary computational study to investigate sensitivity of the solutions with respect to these parameter values.  Tests were carried out using a small `test' problem with a time interval $\mathcal T$ of length 2 hours and a set $\mathcal F$ consisting of 60 flights with randomly-generated pre-scheduled arrival times and weight classes.  Based on the results of this study, we selected the following parameter values: \\
\begin{itemize}
\item (*) Default population size: $S=20$
\item (*) Default sequence length: $l=15$
\item (*) Minimum population size: $S_{\text{min}}=10$
\item (*) Threshold for performing mutations: $m_{\text{mut}}=25$
\item Threshold for performing ranking and selection: $r=50$
\item Threshold for releasing aircraft from pool: $n_{\text{rel}}=50$
\item Threshold for filtering population: $n_{\text{repop}}=500$
\item Maximum `idle runway probability' for delaying pool release: $\lambda=0$
\item Step sizes for sequence cost estimates: $\psi_n=1/n$\\
\end{itemize}

 \text{ }\indent The parameters marked (*) are used by both the SimHeur and DetHeur algorithms.  The others are used only by the SimHeur algorithm.

\end{document}